\PassOptionsToPackage{table}{xcolor}
\documentclass[leqno,oneeqnum,onefignum,onetabnum,final]{siamart1116}

\title{PDE-Constrained Optimization in Medical Image Analysis\thanks{Funding: This material is based upon work supported by AFOSR grants FA9550-17-1-0190; by NSF grant CCF-1337393; by the U.S. Department of Energy, Office of Science, Office of Advanced Scientific Computing Research, Applied Mathematics program under Award Numbers DE-SC0010518 and DE-SC0009286; by NIH grant 10042242; by DARPA grant W911NF-115-2-0121; and by the Technische Universit\"{a}t M\"{u}nchen---Institute for Advanced Study, funded by the German Excellence Initiative (and the European Union Seventh Framework Programme under grant agreement 291763). Any opinions, findings, and conclusions or recommendations expressed herein are those of the authors and do not necessarily reflect the views of the AFOSR, DOE, NIH, DARPA, and NSF. Computing time on the Texas Advanced Computing Centers Stampede system was provided by an allocation from TACC and the NSF.}}

\author{Andreas Mang\thanks{University of Houston, Department of Mathematics, Houston, TX 77204-3008, US, \href{mailto:andreas@math.uh.edu}{andreas@math.uh.edu}}
\and Amir Gholami\thanks{Department of Electrical Engineering and Computer Sciences at University of California, Berkeley, CA 94720-1770, \href{mailto:amir@accfft.org}{amir@accfft.org}}
\and Christos Davatzikos\thanks{Department of Radiology, University of Pennsylvania, Philadelphia, PA 19104-2643, US, \href{mailto:christos.davatzikos@uphs.upenn.edu}{christos.davatzikos@uphs.upenn.edu}}
\and George Biros\thanks{Institute for Computational Engineering and Sciences, University of Texas at Austin, Austin, TX 78712-0027, US, \href{mailto:gbiros@acm.org}{gbiros@acm.org}}}

\usepackage{amsmath,amssymb,mathtools}
\usepackage[sc]{mathpazo}
\usepackage{algorithm,algorithmic}
\usepackage[sort,compress]{cite}
\usepackage{graphicx}
\usepackage{multirow}
\usepackage{makecell}
\usepackage{textcomp}
\usepackage{rotating}
\usepackage{siunitx}
\usepackage{paralist}
\usepackage{booktabs}
\usepackage{lscape}
\usepackage{geometry}
\usepackage[hyperpageref]{backref}
\usepackage[table]{xcolor}

\newtheorem{remark}{Remark}

\definecolor{red}{RGB} {200, 16, 46} 
\makeatletter
\@ifclassloaded{standalone}{}{
\hypersetup{colorlinks,
pdfauthor={Andreas Mang et al.},
linkcolor=red,urlcolor=red,
citecolor=red}}
\makeatother%

\newcommand{\figref}[1]{Fig.~\ref{#1}}
\newcommand{\tabref}[1]{Tab.~\ref{#1}}
\newcommand{\secref}[1]{\S\ref{#1}}
\newcommand{\algref}[1]{Alg.~\ref{#1}}
\newcommand{\remref}[1]{Rem.~\ref{#1}}

\newcounter{runidnum}

\newcolumntype{R}{>{\columncolor{gray!20}}r}
\newcolumntype{L}{>{\columncolor{gray!20}}l}
\newcolumntype{C}{>{\columncolor{gray!20}}c}

\newcommand{\algadjust}{\centering\small\renewcommand\arraystretch{1.2}}

\newcommand{\vect}[1]{\ensuremath{\boldsymbol{#1}}}             

\newcommand{\F}[1]{\ensuremath{\mathcal{#1}}}                   
\newcommand{\D}[1]{\ensuremath{\mathcal{#1}}}                   
\newcommand{\ns}[1]{\ensuremath{\mathbf{#1}}}                   
\newcommand{\fs}[1]{\ensuremath{\mathcal{#1}}}                  
\newcommand{\idiv}{\ensuremath{\nabla\cdot}}                    
\newcommand{\igrad}{\ensuremath{\nabla}}                        
\newcommand{\ilap}{\rotatebox[origin=c]{180}{$\nabla$}}         

\newcommand{\half}[1]{\frac{#1}{2}}
\renewcommand{\d}[1]{\mathop{}\!\mathrm{d}#1}

\newcommand{\p} {\partial}

\newcommand{\defeq}{\ensuremath{\mathrel{\mathop:}=}}

\newcommand{\T}{\ensuremath{\mathsf{T}}}

\newcommand{\di}[1]{\ensuremath{\mathbf{#1}}}                    
\newcommand{\digr}[1]{\ensuremath{\pmb{#1}}}

\newcommand{\bipa}{\begin{inparaenum}[(\itshape i\upshape)]}
\newcommand{\eipa}{\end{inparaenum}}
\newcommand{\bipasub}{\begin{inparaenum}[(\itshape a\upshape)]}
\newcommand{\eipasub}{\end{inparaenum}}

\newcommand{\ipoint}[1]{\textit{\textbf{#1}}}

\reversemarginpar
\newcommand{\mmargin}[1]{{\marginpar{\em\tiny #1}}}\renewcommand{\mmargin}[1]{}

\graphicspath{{figures/}}

\begin{document}

\maketitle

\begin{abstract}
PDE-constrained optimization problems find many applications in medical image analysis, for example, neuroimaging, cardiovascular imaging, and oncological imaging. We review related literature and give examples on the formulation, discretization, and numerical solution of PDE-constrained optimization problems for medical imaging. We discuss three examples. The first one is image registration. The second one is data assimilation for brain tumor patients, and the third one data assimilation in cardiovascular imaging. The image registration problem is a classical task in medical image analysis and seeks to find pointwise correspondences between two or more images. The data assimilation problems use a PDE-constrained formulation to link a biophysical model to patient-specific data obtained from medical images. The associated optimality systems turn out to be sets of nonlinear, multicomponent PDEs that are challenging to solve in an efficient way.

The ultimate goal of our work is the design of inversion methods that integrate complementary data, and rigorously follow mathematical and physical principles, in an attempt to support clinical decision making. This requires reliable, high-fidelity algorithms with a short time-to-solution. This task is complicated by model and data uncertainties, and by the fact that PDE-constrained optimization problems are ill-posed in nature, and in general yield high-dimensional, severely ill-conditioned systems after discretization. These features make regularization, effective preconditioners, and iterative solvers that, in many cases, have to be implemented on distributed-memory architectures to be practical, a prerequisite. We showcase state-of-the-art techniques in scientific computing to tackle these challenges.
\end{abstract}

\newcommand{\slugmaster}{\slugger{siims}{xxxx}{xx}{x}{x--x}}


\pagestyle{myheadings}
\thispagestyle{plain}
\markboth
{A. MANG, A. GHOLAMI, C. DAVATZIKOS, AND G. BIROS}
{PDE-CONSTRAINED OPTIMIZATION IN MEDICAL IMAGE ANALYSIS}

\section{Introduction}
\label{s:intro}

We review methods that integrate medical image analysis algorithms with biophysical models. There are several benefits in integrating biophysical modeling with medical imaging. It can aid prognosis, diagnosis, and the design of new treatment protocols. It can aid our understanding of how biophysics relate to imaging information. Examples of medical image analysis algorithms include image segmentation, image registration, and parameter estimation for biomarkers.

The literature on the integration of biophysical (or biophysically inspired) models with image analysis algorithms is quite extensive. We will consider three problems: The first one is diffeomorphic image registration (see \secref{s:diffeomorphic-registration}). This problem is paramount to many applications in medical imaging~\cite{Fischer:2008a,Sotiras:2013a}. It is about establishing a pointwise spatial correspondence $y:\ns{R}^d\rightarrow\ns{R}^d$, $d\in\{2,3\}$, between two images $m_R:\ns{R}^d\rightarrow\ns{R}$ and $m_T:\ns{R}^d\rightarrow\ns{R}$ of the same object so that the transformed \emph{template image} $m_T\circ y$ becomes similar to the \emph{reference image} $m_R$, i.e,. $m_R(x) \approx (m_T\circ y)(x)$ for all $x$, where $\circ$ denotes the function composition~\cite{Fischer:2008a,Modersitzki:2009a}. We require that the map $y$ is a \emph{diffeomorphism}, i.e., $y$ is a bijection, continuously differentiable, and has a continuously differentiable inverse. In our formulation, we do not directly invert for the map $y$ but introduce a pseudo-time variable $t\in[0,1]$ and invert for the velocity $v : \ns{R}^d \times [0,1] \rightarrow\ns{R}^d$ of $y$. In the context of Lagrangian methods, the map $y(x)$ at a particular location $x$ represents the end point of the characteristic at $t=1$ defined by the velocity $v(x,t)$ (see~\cite{Mang:2017a}). In the present work, we use an Eulerian formulation instead. Here, the PDE constraint is, in its simplest form, given by a hyperbolic transport equation for the image intensities $m_T(x)$ subjected to the velocity $v(x,t)$. The regularization model is a Sobolev norm that stipulates smoothness requirements on $v$. If $v$ is sufficiently smooth it is guaranteed that $v$ gives rise to a diffeomorphism $y$~\cite{Beg:2005a,Chen:2011a,Younes:2010a}. This formulation can be augmented by biophysics operators (hard constraints\footnote{In general, constraints can either be hard or soft. Hard constraints are a set of conditions that the variables are required to satisfy. Soft constraints are penalties for variables that appear in the objective functional; they penalize a deviation of the variables from a condition.}) to incorporate additional prior knowledge on the expected deformation map. One motivation for adding such constraints is that often, especially in longitudinal studies of the same patient, there is a real deformation by the realization of an actual physical phenomenon. A simple constraint, which nonetheless poses significant numerical challenges, is the incompressibility of tissue~\cite{Mang:2015a,Mang:2016a,Mang:2016c}. Examples for more complex biophysical constraints are brain tumor growth models~\cite{Hogea:2008a,Gooya:2012a,Scheufele:2017a,Zacharaki:2008b,Zacharaki:2008a,Zacharaki:2009a} or cardiac motion models~\cite{Sundar:2009a}.

We will discuss all of these constraints either in the context of diffeomorphic image registration or in the more generic context of data assimilation in medical imaging. This brings us to the second problem we are addressing in this article: Data assimilation in brain tumor imaging~\cite{Gholami:2016a,Hogea:2008b,Konukoglu:2010b,Mang:2012b,Menze:2011a} (see \secref{s:dataassim-tumor}). The PDE constraint is, in its simplest form, a nonlinear parabolic differential equation. The inversion variables are, e.g., the initial condition, the growth rate of the tumor, or a diffusion coefficient that controls the net migration of cancerous cells within brain parenchyma. The regularization model is in our case an $L^2$-penalty. The third problem is cardiac motion estimation (see \secref{s:heart}).

\subsection{Contributions}

\begin{itemize}
\item We provide a comprehensive review of existing work on three topics for PDE-constrained optimization with application to medical image analysis:
  \begin{enumerate}
  \item diffeomorphic image registration (see \secref{s:diffeomorphic-registration}),
  \item data assimilation in brain tumor imaging (see \secref{s:dataassim-tumor}), and
  \item data assimilation in cardiac imaging (see \secref{s:heart}).
  \end{enumerate}
\item We discuss the implementation of a distributed-memory Newton--Krylov solver for two of these problems (see \secref{s:numerics}).
\end{itemize}

In this paper, we summarize our recent contributions in the field of diffeomorphic image registration~\cite{Mang:2015a,Mang:2016a,Mang:2016c,Mang:2017a,Mang:2017c} and data assimilation in brain tumor imaging~\cite{Gholami:2016a,Gholami:2017a,Mang:2012b,Mang:2014a}. In particular, we will present the formulations described in~\cite{Mang:2015a,Mang:2016a} and~\cite{Gholami:2016a}, the Newton--Krylov methods presented in~\cite{Mang:2015a,Mang:2017c} and~\cite{Gholami:2016a}, and review the distributed-memory implementations described in~\cite{Mang:2016a,Gholami:2017a}. We have also added some new results, which have not been presented elsewhere (see \secref{s:experiments:reg}).

\subsection{Mathematical Setting}

The problems we consider in this paper can be viewed as an inverse problem of estimating a parameter function $w$ (\emph{inversion field}) based on observations $b$, which we assume to be solutions of a system of PDEs $\D{C}[m,w] = 0$ with $\D{C}[m,w]\defeq \D{A}(w)m - q$, defined on a domain $\bar{\Omega}\times[0,T]$, $\Omega=[0,1]^d$, $d\in\{2,3\}$. We, e.g., consider linear hyperbolic and nonlinear parabolic PDE constraints $\D{C}$. We assume that the data $b$ is a nonlinear function of the parameters $w$, i.e., the parameter-to-observation map is given by $b = \D{Q}\D{A}(w)^{-1}q + \delta$, where $\delta$  is measurement noise. Here, $\D{Q}$ is an observation operator, i.e., a projection of the \emph{state variable} $m$ onto locations in $\Omega$ (or $\Omega \times [0,T]$) to which the observation (or desired state) $b$ is associated. The problem of recovering $w$ from $b$ is often ill-posed due to the fact that the measured data is finite and perturbed by noise. A common approach to tackle this challenge is to solve for a smooth, locally unique solution of a nearby problem by imposing prior knowledge based on an adequate regularization scheme~\cite{Engl:1996a}, for example a Tikhonov functional $\D{R}$. Overall, we arrive at the following formulation
\begin{align}
  \label{e:pde-opt}
  \min_{m,w}\D{J}[m,w] &= \half{1}\|\D{Q}m - b\|^2_{L^2(\Omega)} + \beta\,\D{R}[w]
  \\
  \text{subject to}\;\; \D{C}[m,w] &= 0.\nonumber
\end{align}

The problem in~\eqref{e:pde-opt} balances data fidelity (for simplicity, we only consider an $L^2$-distance) with regularity (controlled by the (regularization) parameter $\beta > 0$). Designing efficient algorithms for PDE-constrained optimization problems remains a significant challenge~\cite{Akcelik:2006a,Biegler:2003a,Borzi:2012a,Herzog:2010a,Hinze:2009a,Leugering:2014a}. These problems typically involve an infinite number of unknowns, which, upon discretization, leads to high-dimensional systems (in our case, e.g., $\D{O}(\num{1e+07})$ unknowns for clinically relevant problems). Usually, the solver has to be tailored to the structure of the operators. We focus on derivative/adjoint based algorithms for PDE-constrained optimization. The evaluation of the objective functional and its derivatives involves expensive computations, which, in our examples, are the solution of PDEs. We tackle these challenges by designing computational methods that rigorously follow mathematical and physical principles, and are based on efficient, robust, high-fidelity algorithms with a sound theoretical basis.

\begin{remark}\label{r:bayesian-inversion}
Typically, the model parameters $w$ are \emph{not} the final quantity of interest. More often, we need to predict some other future quantity (e.g., tumor infiltration to healthy tissue or the probability of tumor recurrence). Using a deterministic formulation, like the one in~\eqref{e:pde-opt}, only provides point estimates. Due to uncertainties in the data $b$ (denoted above by $\delta$), the model parameters $w$, the inversion algorithm, and the mathematical model $\fs{C}$, we require confidence intervals for the quantities of interest, not just point estimates. That is, we are interested in uncertainty propagation from the input to the quantity of interest. This can be achieved by using a Bayesian framework~\cite{Kaipio:2005a,Tarantola:2005a,Sullivan:2015a}. Formulating the problem as a Bayesian inverse problem provides the posterior joint probability density of the parameters $w$. Related quantities of interest can be computed by computing expectations with respect to the posterior. The solution of the deterministic inverse problem corresponds to the maximum a posteriori probability estimate. This connection between deterministic inversion and statistical inference can be used to devise efficient sampling strategies that exploit local problem structure~\cite{Geweke:1999a,Geweke:2003a,Martin:2012a,Petra:2014a}. In what follows, we will focus on deterministic inversion.
\end{remark}

\subsection{Limitations and Open Issues}

Despite the success of biomechanical models, in many cases their clinical application is limited to qualitative analysis. Biomechanical simulators are hampered by imprecise and complex constitutive laws, and the uncertainties in boundary and initial conditions. Moreover, even assuming we had a complete knowledge of the models and their parameters, solving the underlying PDEs is computationally challenging; accounting for the uncertainty amplifies the computational difficulties.

The key challenges and limitations of biomechanically informed algorithms for medical image analysis can be summarized as follows:
\begin{itemize}
\item Models are often phenomenological; many parameters in the equations do not have direct biological counterparts and can therefore not be directly obtained from data.
\item Boundary and initial conditions have significant uncertainties.
\item In a clinical setting, in-vivo medical imaging data has typically poor resolution and may contain significant amount of noise, which makes the calibration of complicated (multi-parameter) models quite challenging.\footnote{We note that high-resolution imaging technologies have emerged; an example is CLARITY imaging~\cite{Chung:2013a,Tomer:2014a,Kutten:2017a}; we revisit this ex-situ imaging technique briefly in~\secref{s:diffeomorphic-registration}. High-resolution imaging techniques are not yet available in a clinical setting. Resolution levels for routinely collected imaging is between \SI{0.5}{mm} and \SI{5}{mm} along each spatial direction (depending on the imaging modality).}
\item The imaging operator is extremely complicated. The question on how observations in imaging data map to model outputs is a difficult one and depends on the scanner. Oftentimes, pseudo-correspondences are used.
\item Complicated models result in ill-conditioned, multi-physics operators that are challenging to solve.
\item We obtain complex optimization problems that are highly nonlinear and time-dependent.
\end{itemize}

We limit this review to problems we have considered in our past work. We note that PDE-constrained optimization appears in other contexts in medical imaging. For example, many image reconstruction algorithms involve physical models and the reconstruction phase could be formulated as a PDE-constrained optimization problem. However, this usually leads to overly complicated algorithms. Typically approximation and spectral analysis of the underlying operators are preferable and lead to direct reconstruction algorithms. Here, we briefly give some references for problems in which PDE-constrained optimization has been considered in practice. In elastography, the goal is to estimate elastic material properties using, e.g., ultrasound. The main PDE is wave propagation. Examples can be found in~\cite{Goenezen:2012a,Oberai:2003a,Ophir:1999a}. In optical diffusion tomography, the goal is to image optical properties of tissue using diffuse approximations of light scattering. The main PDEs are diffusion, radiation-diffusion, or highly-damped Helmholtz. Examples include~\cite{Arridge:1999a,Arridge:2009a,Joshi:2004a,Ren:2006a,Saratoon:2013a}.

\subsection{Organization and Notation}

We denote the standard $L^2$ inner product by $\langle\cdot,\cdot\rangle_{L^2(\Omega)}$ and the $L^2$ norm by $\|\cdot\|_{L^2(\Omega)}$, both defined on the spatial domain $\Omega\subset\ns{R}^d$, $d\in\{2,3\}$. We denote discretized quantities by upright boldface letters. That is, the discrete representation of $y$ is denoted by $\di{y}$. Likewise, a superscript $h$ is added to the operators whenever we refer to discretized quantities.

We review the literature and present the formulation for the diffeomorphic registration problem in \secref{s:diffeomorphic-registration} and for the data assimilation problem in brain tumor imaging in \secref{s:dataassim-tumor}. We review relevant literature for cardiac motion estimation in \secref{s:heart}. We present the implementation of our solver and its deployment in high-performance computing platforms in \secref{s:numerics}. We showcase exemplary numerical results in~\secref{s:numerical-experiments}. We conclude with \secref{s:conclusions}. Additional details about our solver can be found in the appendix (\secref{s:gn-krylov-method} and \secref{s:newton-step}).

\section{Diffeomorphic Registration}
\label{s:diffeomorphic-registration}

Let us start by defining more precisely what we mean by diffeomorphic image registration. In diffeomorphic image registration we require that the spatial transformation $y$ that maps the template image $m_T$ to the reference image $m_R$ is a \emph{diffeomorphism}. A diffeomorphism is an invertible map, which is continuously differentiable (in particular a $C^1$-function) and maps $\Omega$ onto itself. Formally, we require that $\det \igrad y$ does not change sign or vanish for $y$ to be locally diffeomorphic.

\subsection{Literature Review}
\label{s:literature-review-reg}

Providing a comprehensive review on image registration is impossible due to the large body of literature. We refer to~\cite{Fischer:2008a,Modersitzki:2004a,Modersitzki:2009a,Sotiras:2013a} for a general introduction to the problem, recent algorithmic developments, and applications. We consider the problem of \emph{diffeomorphic} image registration as a \emph{hyperbolic optimal control problem}~\cite{Borzi:2002a,Chen:2011a,Hart:2009a,Mang:2015a,Mang:2016a,Mang:2017c,Mang:2016c}; i.e., problems with a transport or continuity equation as PDE constraint (we refer to \secref{s:formulation-reg} for details on the problem formulation).\footnote{Similar formulations can be found in other applications, e.g., geophysical sciences~\cite{Fohring:2014a}.} These types of formulations are based on conceptual ideas from computational fluid dynamics and flow control~\cite{Gunzburger:2003a}. They share characteristics with traditional optical flow formulations~\cite{Horn:1981a,Kalmoun:2011a,Pock:2007a,Ruhnau:2007a} and formulations for optimal mass transport~\cite{Angenent:2003a,Benamou:2000a,Benzi:2011a,Haber:2015a,Simoncini:2012a,Rehman:2009a}.

Viscous fluid formulations for diffeomorphic image registration are pioneered in~\cite{Christensen:1994a,Christensen:1996a}. Here, the basic idea is to introduce a pseudo-time variable and parametrize the deformation map $y$ by its velocity $v$. This approach has been embedded into a variational framework in~\cite{Beg:2005a,Dupuis:1998a,Miller:2004a,Miller:2001a,Trouve:1998a}; the problem formulation takes the form
\[
\min_v \half{1}\|m_T\circ \phi_1^{-1}-m_R\|^2_{L^2(\Omega)}+\half{\beta}\int_0^1\|v\|^2_{\fs{S}}\d t
\]

\noindent subject to $\d_t \phi = v(\phi)$ for $t\in(0,1]$ and $\phi = \operatorname{id}$ for $t=0$, where $\phi_1 \defeq \phi(t=1)$ is equivalent to $y^{-1}$, and $\operatorname{id}(x) = x$. In the formulation above, we seek to find a time dependent, smooth, and compactly supported velocity field $v\in L^2([0,1],\fs{S})$, where $\fs{S}$ is a Sobolev space of certain regularity. The Sobolev norm in the problem formulation above guarantees that $v$ belongs to a space $\fs{S}$ of regular vector fields, which in turn guarantees (the regularity requirements are discussed in \cite{Dupuis:1998a,Trouve:1998a,Beg:2005a}) that the solution to $\d_t\phi = v(\phi)$ is a diffeomorphism. This approach embeds the diffeomorphisms in a Riemannian space; the time integral of the squared Sobolev norm of $v$ is a geodesic distance (a Riemannian metric) between the identity map $\operatorname{id}$ at $t=0$ and the diffeomorphism $\phi_1$ at $t=1$. The optimal $v$ defines the shortest path on the manifold of diffeomorphisms that connects $\operatorname{id}$ with $\phi_1$.

As we will see below, we will use a PDE-constrained (Eulerian) formulation instead. A related formulation that includes additional hard and soft constraints on $v$ and $\det \igrad y$ can be found in~\cite{Lee:2010a,Lee:2011a}. This is different to the PDE-constrained formulations in~\cite{Borzi:2002a,Chen:2011b,Chen:2011a,Hart:2009a,Mang:2015a,Mang:2016a,Mang:2017c,Mang:2016c}. Here, a linear advection equation is used to model the transport of the intensities of $m_T$. The latter approach does not operate on the space of diffeomorphisms; $y$ does not appear explicitly, only $v$ does. This formulation has been augmented by hard or soft constraints on the divergence of $v$~\cite{Borzi:2002a,Chen:2011a,Mang:2015a,Mang:2016a,Mang:2016c,Ruhnau:2007a}, rendering the deformation map $y$ (nearly) incompressible~\cite[p.~77ff.]{Gurtin:1981a}.\footnote{An interesting direction of research is to augment these PDE constraints by more complex biophysics operators~\cite{Hogea:2008a,Gooya:2012a,Sundar:2009a,Zacharaki:2008b,Zacharaki:2008a,Zacharaki:2009a}. This results in additional parameters that need to be calibrated.} In both approaches $v$ has to be sufficiently smooth to guarantee that the associated deformation map $y$ is a diffeomorphism. These smoothness requirements are typically stipulated by the regularization operator---a Sobolev norm, the order of which depends on the smoothness of the images, the particular form of the transport equation, and the existence of additional constraints (such as, e.g., a divergence-free velocity)~\cite{Barbu:2016a,Beg:2005a,Chen:2011b,Chen:2011a,Crippa:2007a,Dupuis:1998a,Trouve:1998a} (see~\secref{s:formulation-reg}).

The differences between our formulation and PDE-constrained formulations for optimal mass transport~\cite{Benamou:2000a,Benzi:2011a,Haber:2015a,Simoncini:2012a} are the constraint (a continuity equation to preserve mass\footnote{Applications for mass-preserving registration can be found in~\cite{Burger:2013a,Mang:2017a,Wlazlo:2016a}.} versus a transport equation) and the regularization model ($L^2$-norm of $v$ scaled by the transported density versus Sobolev-norm of $v$). The conceptual differences between our formulation and traditional optical flow formulations~\cite{Horn:1981a,Kalmoun:2011a,Pock:2007a,Ruhnau:2007a} are that we do not consider a time series of images and the transport equation is a hard constraint; it does not directly enter the objective functional. PDE-constrained formulations for optical flow, which are equivalent to our formulation, are described in~\cite{Andreev:2015a,Barbu:2016a,Borzi:2002a,Chen:2011a}.

There has only been a limited amount of work devoted to the design of effective algorithms for velocity-based diffeomorphic image registration that exploit state-of-the-art technology in scientific computing. There are several challenges for the design of efficient solvers. First and foremost, image registration is an ill-posed, highly nonlinear inverse problem. It is established that gradient descent schemes for numerical optimization (which are still predominantly used~\cite{Avants:2006a,Avants:2008a,Avants:2011a,Beg:2005a,Chen:2011a,Ha:2010a,Hart:2009a,Vialard:2012a}) only converge slowly when applied to nonlinear, ill-posed problems. Recently, several groups have developed Newton-type algorithms to address this issue~\cite{Ashburner:2011a,Benzi:2011a,Mang:2015a,Mang:2016a,Mang:2017c,Mang:2016c,Mang:2017a,Simoncini:2012a,Herzog:2018a}. In this context, the design of an effective preconditioner for the Hessian is essential for these methods to be competitive in terms of the time-to-solution; different approaches have been discussed in~\cite{Benzi:2011a,Mang:2015a,Mang:2017c,Mang:2017a,Simoncini:2012a,Herzog:2018a} (see~\cite{Benzi:2005a} for an introduction to preconditioning of saddle point problems).
Another issue is the size of the search space, in particular when considering 3D problems (in space).\footnote{The number of unknowns is the number of discretization points in space times the number of discretization points in time for the velocity times the dimensionality of the ambient space (we invert for a time-dependent vector field).} The size of the search space can be reduced by exploiting the fact that the momentum of the variational problem is constant in time (in Lagrangian coordinates). That is, the flow of the diffeomorphism can completely be encoded by the momentum at $t=0$~\cite{Ashburner:2011a,Miller:2006a,Vialard:2012a,Younes:2007a,Younes:2010a}. Another approach is to invert for a stationary velocity field $v(x)$~\cite{Arsigny:2006a,Ashburner:2007a,Hernandez:2009a,Mang:2016a,Mang:2016c,Mang:2017a}.\footnote{A stationary velocity field $v$ is a velocity field that is constant in time, as opposed to a nonstationary---i.e., time-dependent or transient---velocity field. We note that stationary velocity fields do not cover the entire space of diffeomorphisms and do not provide a Riemannian metric on this space, something that may be desirable in certain applications~\cite{Beg:2005a,Miller:2004a,Zhang:2015a}; this requires time-dependent velocities. The work in~\cite{Mang:2015a} uses a Galerkin method to control the number of unknowns in time. It is shown experimentally that stationary and nonstationary velocities yield an equivalent registration quality in terms of data mismatch.}
In addition to that, we have to deal with an expensive parameter-to-observation map, which, in our formulation, is a hyperbolic transport equation. Several strategies to efficiently solve these equations have been considered; they range from conditionally stable total variation diminishing~\cite{Borzi:2002a}, second order Runge-Kutta~\cite{Mang:2015a,Mang:2016a}, and high-order essentially nonoscillatory~\cite{Hart:2009a} schemes, to unconditionally stable implicit Lax-Friedrich~\cite{Benzi:2011a,Simoncini:2012a}, semi-Lagrangian~\cite{Beg:2005a,Chen:2011a,Mang:2017c,Mang:2016c}, and Lagrangian~\cite{Mang:2017a} schemes.

We note, that there exists a large body of literature that discusses effective parallel implementations of solvers for PDE-constrained optimization problems; examples include~\cite{Akcelik:2002a,Akcelik:2006a,Biros:1999a,Biros:2005a,Biros:2005b,Biegler:2003a,Biegler:2007a,Schenk:2009a}. There has been only little work on distributed-memory algorithms for diffeomorphic image registration. Real-time implementations for low dimensional (rigid) registration are, e.g., described in~\cite{Shams:2010b,Ruehaak:2017a}. Popular software packages for (in many cases) diffeomorphic image registration are described in~\cite{Ashburner:2007a,Avants:2011a,Klein:2010a,Modat:2010a,Modersitzki:2009a,Ou:2011a,Vercauteren:2008a,Vercauteren:2009a}. None of these uses a PDE-constrained formulation like the one discussed here. Most implementations use open multi-processing~\cite{Avants:2011a,Klein:2010a,Vercauteren:2008a,Vercauteren:2009a} or graphic processing units (GPUs) for acceleration~\cite{Gu:2010a,Ha:2009a,Modat:2010a,Rehman:2009a,Sommer:2011a,MuyanOzcelik:2008a,Gu:2010a,Ellingwood:2016a,Ko:2008a,ValeroLara:2014a,Luo:2015a}. A survey of registration algorithms that exploit GPUs for acceleration can be found in~\cite{Fluck:2011a}. A literature survey on high performance computing in image registration can be found in~\cite{Shams:2010a,Eklund:2013a,Shackleford:2013a}. Works on deploying solvers for diffeomorphic image registration to supercomputing platforms can be found in~\cite{Ha:2009a,Ha:2010a,Ha:2011a,Mang:2016c,Gholami:2017a}. The work that is most closely related to ours~\cite{Mang:2016c,Gholami:2017a} is~\cite{Ha:2009a,Ha:2011a} and~\cite{ValeroLara:2014a}. In~\cite{Ha:2009a,Ha:2010a,Ha:2011a} a multi-GPU accelerated implementation of the diffeomorphic image registration approach described in~\cite{Joshi:2005a} is presented. The work in~\cite{ValeroLara:2014a} discusses a GPU accelerated implementation of the diffeomorphic algorithm described in~\cite{Ashburner:2007a}.

Recent advances in medical imaging~\cite{Chung:2013a,Tomer:2014a,Kutten:2017a} result in data sizes that pose the demand for deploying scalable, distributed-memory implementations.\footnote{The resolution of these datasets can reach $\SI{5}{\micro\metre}\times\SI{5}{\micro\metre}\times\SI{5}{\micro\metre}$ resulting in $\mathcal{O}(\SI{4.8}{\tera\byte})$ of data (if stored with half precision)~\cite{Tomer:2014a}. Overall, this results in 2.4 trillion discretization points in space.} To the best of our knowledge, only our solver~\cite{Mang:2016c,Gholami:2017a} has been shown to scale up to similar problem sizes. Moreover, we are the first group that proposed an implementation of a globalized (Gauss--)Newton--Krylov solver for a PDE-constrained approach to diffeomorphic image registration that scales on supercomputing platforms~\cite{Mang:2016c,Gholami:2017a}.

\begin{figure}
\centering
\includegraphics[width=0.9\textwidth]{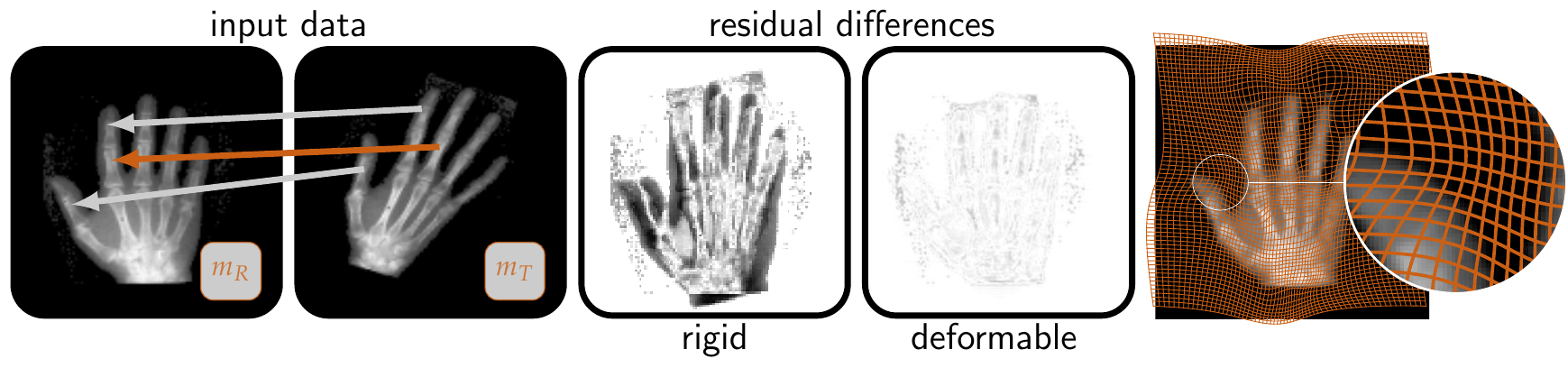}
\caption{Diffeomorphic image registration problem. The inputs to this inverse problem are scalar intensity values of two images---the \emph{reference image} $m_R(x)$ and the \emph{template image} $m_T(x)$ (image to be registered)---of the same object (left); data taken from~\cite{Amit:1994a,Modersitzki:2009a}. The inherent assumption is that there exists a geometric transformation (e.g., a rigid body transformation (middle) or a deformable transformation (middle right and right)) that relates points in $m_R$ to points in $m_T$ (arrows). These spatial correspondences are, in our formulation, solely based on intensity values; points with similar intensities map to one another. In general, this can result in mappings that are not meaningful (orange arrow). In our formulation, we consider a mapping \emph{meaningful} or \emph{plausible} if it is a diffeomorphism, i.e., it is smooth, onto, one-to-one, and it has a smooth inverse.\label{f:correspondence-problem}}
\end{figure}

\subsection{Formulation}
\label{s:formulation-reg}

We assume that the reference image $m_R : \bar{\Omega} \rightarrow [0,1]$ and the template image $m_T : \bar{\Omega} \rightarrow [0,1]$ are continuously differentiable functions, and compactly supported on $\bar{\Omega} \defeq (0,2\pi]^d$, $d\in\{2,3\}$. In diffeomorphic image registration we seek to find a diffeomorphism $y$ such that the transformed template image $m_T$ becomes similar to the reference image $m_R$, i.e., $\forall x \in \Omega: (m_T \circ y)(x) \approx m_R(x)$~\cite{Modersitzki:2004a,Modersitzki:2009a}. In our formulation, we do not directly invert for the deformation map $y$. We introduce a pseudo-time variable $t\in[0,1]$ and invert for a smooth velocity field (control variable) $v \in L^2([0,1], \fs{S})$, $\fs{S}\subseteq H^k(\Omega)$, $k\in\{1,2,3\}$, where $H^k$ is a Sobolev space of order $k$ (see below); the velocity $v$ parameterizes the deformation map $y$. We use an Eulerian formulation; in its simplest form, the PDE constraint $\D{C}$ in~\eqref{e:pde-opt}, i.e., the forward operator, is given by the transport equation~\cite{Mang:2015a,Hart:2009a}
\begin{subequations}\label{e:reg:fwd}
\begin{align}
\p_t m + v \cdot \igrad m & = 0   && {\rm in}\;\Omega\times(0,1],\label{e:reg:fwd-pde}\\
                        m & = m_T && {\rm in}\;\Omega\times\{0\},\label{e:reg:fwd-initial}
\end{align}
\end{subequations}

\noindent and periodic boundary conditions on $\p\Omega$. In our formulation, the map $y$ does not appear explicitly. This is different from map-based formulations (i.e., formulations that directly invert for $y$)~\cite{Burger:2013a,Fischer:2008a,Mang:2012a,Modersitzki:2009a} and Lagrangian formulations for velocity-based diffeomorphic image registration~\cite{Avants:2011a,Beg:2005a,Mang:2017a}. In our formulation, we model the deformed template image as the solution of the transport equation~\eqref{e:reg:fwd} at $t=1$ and denote it by $m_1(x) \defeq m(x,t = 1)$.

In the inverse problem we seek to find a velocity $v$ such that the transported template image $m_1(x)$ is equal to $m_R(x)$ for all $x\in\Omega$. As can be seen in~\eqref{e:pde-opt}, we use a squared $L^2$-distance to measure the discrepancy between $m_1(x)$ and $m_R(x)$, where $m_R$ corresponds to $b$ and $Q$ projects the state $m$ to the terminal time frame at $t = 1$, i.e., $m_1 = Qm$.\footnote{Other distance measures, such as mutual information, normalized gradient fields, or cross correlation can be used (see \cite{Modersitzki:2004a,Modersitzki:2009a} for an overview). In our formulation, changing the distance measure will affect the first term in the objective functional and the terminal condition of the adjoint equation~\eqref{e:reg:adjoint} in~\secref{s:optimality-systems:reg}.} Notice, that the reference image $m_R(x)$ as well as the initial condition of the forward problem in \eqref{e:reg:fwd} are measured function, and, hence, may contain noise.

This formulation has been augmented by constraints on the divergence of $v$~\cite{Chen:2011a,Mang:2015a,Mang:2017c,Mang:2016c}. We obtain an incompressible diffeomorphism $y$ if $v$ is divergence-free~\cite[p.~77ff.]{Gurtin:1981a}, i.e., $\det\igrad y = 1$. In this case, $\fs{S}$ is the space of divergence-free velocities with Sobolev regularity of order $k$. We can relax the incompressibility assumption by introducing a mass-source $q$, i.e., $\idiv v = q$~\cite{Mang:2016a}. This is equivalent to adding a penalty on the divergence of $v$ to the variational problem~\cite{Borzi:2002a}. In~\cite{Mang:2017a} it has been suggested to replace~\eqref{e:reg:fwd-pde} with the continuity equation, e.g., used in optimal mass transport formulations~\cite{Benamou:2000a,Benzi:2011a,Haber:2015a,Simoncini:2012a} to model the preservation of mass of $m_T$.

A critical ingredient of the variational problem formulation is the regularization operator $\D{R}$; choosing an adequate Sobolev norm guarantees that the deformation map $y$ associated with $v$ is a diffeomorphism~\cite{Beg:2005a,Dupuis:1998a,Miller:2004a,Miller:2001a,Trouve:1998a}. In general format, $\D{R}$ is given by
\begin{equation}\label{e:reg:regularization}
\D{R}[v] = \half{1} \int_0^1 \|v\|_{\fs{V}}^2\d t
= \half{1} \int_0^1 \langle \D{A} v, v\rangle_{L^2(\Omega)^d} \d t.
\end{equation}

\noindent Here, $\D{A} : \fs{V} \rightarrow \fs{V}^\ast$ is a positive definite, self-adjoint differential operator. In~\cite{Beg:2005a} it is shown that we require (slightly more than) $H^2$-regularity for $v$ if we assume that $m_R$ and $m_T$ are $H^1$-functions. Accordingly, $\D{A} = \D{B}^\dag\D{B}$ with $\D{B}\defeq(-\ilap + \gamma)^\kappa I$, $\kappa \geq 1$, has become a popular choice~\cite{Beg:2005a,Hernandez:2009a,Zhang:2015a}. If we model the images as functions of bounded variation instead, we require that $v$ is an $H^3$-function~\cite{Chen:2011a}. This requirement can be relaxed to $L^2$-integrability if we add a diffusion operator to~\eqref{e:reg:fwd-pde}~\cite{Barbu:2016a}. A detailed study of (near) incompressible $H^1$ flows can be found in~\cite{Chen:2011b,Chen:2011a,Crippa:2007a}. The work in~\cite{Borzi:2002a} considers an additional regularization term for $\p_t v$. In~\cite{Mang:2016a,Mang:2016c} $v$ is modeled as a stationary field; the integral in time in~\eqref{e:reg:regularization} can be dropped. Extensions to nonquadratic $H^1$-regularization models are discussed in~\cite{Mang:2016a,Pock:2007a}. In particular, the work in~\cite{Mang:2016a} introduces a nonlinear regularization operator that enables a control (promote or penalize) of shear (the first variation yields a control equation with a viscosity that depends on the strain rate).

\section{Data Assimilation in Brain Tumor Imaging}
\label{s:dataassim-tumor}

\subsection{Literature Review}
\label{s:literature-tumor}

\begin{figure}
\centering
\includegraphics[width=0.9\textwidth]{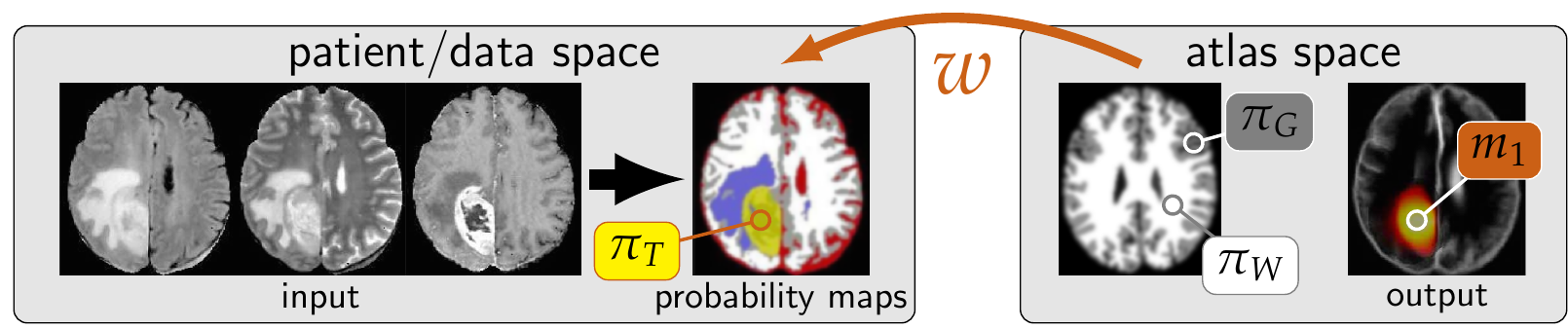}
\caption{Parameter identification/data assimilation in brain tumor imaging. The input to our problem are magnetic resonance images of an individual patient (left; from left to right: fluid attenuation inversion recovery; T2-weighted; and T1-weighted with contrast enhancement). Probability maps for tissue types extracted from these data serve as an input to our inversion; middle; tumor (yellow); edema (purple); white matter (white); gray matter (gray); cerebrospinal fluid (red). We perform the simulations in a standardized atlas space (due to the absence of an estimate of the patient's image without tumor; the data is taken from~\cite{Cocosco:1997a}). This dataset provides probability maps for anatomical regions for a normal brain image (e.g., white and gray matter denoted by $\pi_W(x)$ and $\pi_G(x)$, respectively). The goal is to estimate parameters $w$ so that the predicted state $m_1(x)$ at $t=1$ matches the observed state $\pi_T(x)$ (parameter identification problem). The images are modified from~\cite{Gooya:2012a}.\label{f:data-assimilation-problem}}
\end{figure}

We view biophysics simulations as a powerful tool to augment morphological and functional medical imaging data. There is a long tradition in the design of mathematical models of cancer progression~\cite{Bellomo:2008a,Roose:2007a,Wang:2008a}. In many cases, the models are rather simplistic.\footnote{More sophisticated multi-species models that, e.g., account for hypoxia, necrosis and angiogenesis can be found in~\cite{HawkinsDaarud:2013a,Gu:2012a,Rahman:2017a}.} This is due to the fact that medical imaging only provides scarce information to drive the model (see \figref{f:data-assimilation-problem}); we have to limit the size of the parameter space. Most of the approaches used for data assimilation in medical imaging exclusively rely on observations derived from morphologic imaging data and consider continuum models; cancerous cells are not tracked individually but modeled as a density $m:\bar{\Omega} \times [0,1] \rightarrow [0,1]$ on $\bar{\Omega}\times [0,1] \subset\ns{R}^{d+1}$, $d\in\{2,3\}$. In their simplest form, these models assume that cancer progression is dominated by two phenomena: cell division/mitosis and cell migration, resulting in reaction-diffusion type equations~\cite{Murray:1989a,Swanson:2000a,Swanson:2002a,Jackson:2015a}.\footnote{Models that account for the mechanical interaction of the tumor with its surroundings have been described in~\cite{Chen:2012b,Clatz:2005a,Hogea:2007a,Mohamed:2005a,Weis:2017a,Wong:2015a,Wong:2017a}.} Despite their simplicity and phenomenological character, these types of models have been successfully applied to capture the appearance of tumors in medical images~\cite{Clatz:2005a,Konukoglu:2010a,Liu:2014a,Mang:2014a,Mang:2012b,Rekik:2013a}. They have been used to \bipa\item study tumor growth patterns in individual patients~\cite{Clatz:2005a,Mang:2014a,Mang:2012b,Rekik:2013a,Swanson:2000a}, \item extrapolate the physiological boundary of tumors~\cite{Konukoglu:2010a,Mosayebi:2012a}, \item or study the effects of clinical intervention~\cite{Le:2015a,Le:2017a,Mi:2014a,Powathil:2007a,Rockne:2010a,Swanson:2008a,Weis:2017a}\eipa.

Most of the work mentioned above is concerned with brain tumor imaging; exceptions are, e.g., \cite{Weis:2017a} (breast cancer) and~\cite{Liu:2014a} (pancreatic cancer). A key challenge is the design of methods that have the potential to generate patient-specific predictions. In general, this requires the solution of an inverse parameter identification problem. Several groups have addressed this problem. There has only been little work on model-based image analysis in brain tumor imaging based on PDE-constrained optimization~\cite{Colin:2014a,Gholami:2016a,Hogea:2008a,Knopoff:2013a,Knopoff:2017a,Liu:2014a,Mang:2012b,Mang:2014a,Quiroga:2015a};\footnote{We note that there exists a large body of literature on optimal control problems with similar PDEs as constraints in other areas. Examples can, e.g., be found in~\cite{Barthel:2010a,Croft:2015a,Figueiredo:2011a,Figueiredo:2013a,Pearson:2013a}.} the works in~\cite{Colin:2014a,Gholami:2016a,Hogea:2008a,Knopoff:2013a,Knopoff:2017a,Liu:2014a,Quiroga:2015a} consider adjoint based approaches for numerical optimization. Others use derivative-free optimization~\cite{Chen:2012b,Knopoff:2017a,Konukoglu:2010b,Mang:2012a,Mang:2014a,Mi:2014a,Wong:2015a,Wong:2017a}, finite-difference approximations to the gradient~\cite{Hormuth:2015a}, or tackle the parameter estimation problem within a Bayesian framework~\cite{Le:2015a,Le:2017a,Menze:2011a,HawkinsDaarud:2013b,Oden:2013a,Lima:2016a,Collis:2017a,Lima:2017a} (see \remref{r:bayesian-inversion}). These approaches have in common that they only require the evaluation of an objective functional (i.e., essentially a repeated solution of the forward problem for perturbed parameters). As such, they are easy to implement. However, it is established that they are in general not as efficient as gradient based optimization strategies (at least for high-dimensional problems); they display slow convergence resulting in an excessive number of iterations/evaluations of the forward model. Another strategy that has been considered, is to exploit the asymptotic behavior of reaction-diffusion equations~\cite{Swanson:2008a,Harpold:2007a,Konukoglu:2010b}. The basic idea is to consider the boundary of the abnormality visible in medical imaging as a traveling wave solution and convert it into parameter estimates for the growth/migration rate (see~\cite{Harpold:2007a}). The work in~\cite{Hogea:2008a} is to the best of our knowledge the first to consider adjoint equations in the context of data assimilation in brain tumor imaging. The adjoint equations are derived for the one-dimensional case; they fall back to a derivative-free optimization strategy when solving the three-dimensional problem. The work in~\cite{Colin:2014a} is based on a multi-species model. They derive adjoint equations for inverting for the vascularization. The work in~\cite{Liu:2014a} extends~\cite{Hogea:2008a}. They derive first-order optimality conditions for an advection-reaction-diffusion PDE-constrained optimization problem; there is no information on the numerical treatment of the problem.

In~\cite{Quiroga:2015a} a coupled reaction-diffusion system for the density of normal tissue, cancerous tissue, and ion concentration is considered. They derive adjoint equations for a scalar parameter controlling the excess ion concentration. The work in~\cite{Knopoff:2013a} models the tumor as a radially symmetric spheroid. The only work we are aware of (in the context of model based brain tumor image analysis) that derives Hessian information for numerical optimization is~\cite{Gholami:2016a}. We describe this model in more detail in~\secref{s:formulation-tumor}. Aside from numerical optimization, there has only been little work on the design of effective solvers (in terms of time-to-solution and rate of convergence). Conditionally stable explicit time integrators are a common choice due to their simplicity~\cite{Colin:2014a,Weis:2017a}. Using methods that require a large number of time steps can become computationally prohibitive in a three-dimensional setting not only in terms of the time-to-solution but also due to memory requirements.\footnote{A naive implementation may require storage of the time history of the state and adjoint fields; one remedy is the implementation of check-pointing/domain decomposition schemes~\cite{Akcelik:2002a,Griewank:1992a,Heinkenschloss:2005a}.} The works in~\cite{Hogea:2008a,Gholami:2016a,Mang:2012a,Mang:2014a} consider unconditionally stable schemes.

An open problem in data assimilation in brain tumor imaging is how to establish a \emph{physiologically meaningful} correspondence between the model output (predicted state; typically a cell density) and the data (observed state; typically morphological imaging information, i.e., an abnormality in the images; see \figref{f:data-assimilation-problem} for an example). As a consequence one has to rely on heuristics and pseudo-correspondences. One such example is to relate (manual) segmentations of imaging abnormalities to a detection thresholds for the computed population density of cancerous cells~\cite{Gholami:2016a,Harpold:2007a,Le:2015a,Le:2017a,Mang:2014a,Mang:2012a,Swanson:2008a}. Another approach is to derive cell density estimates from physiological imaging data~\cite{Atuegwu:2012a,Hormuth:2015a,Weis:2017a}. Either of these techniques only provides scarce information to drive the parameter estimation, which in turn limits the size of the search space, and as such the complexity of the models.

Considering parallel, distributed-memory implementations, we are not aware of any work in this area other than the one presented by our group~\cite{Gholami:2017a}.\footnote{We note that there certainly exist parallel implementations of solvers for similar PDE-constrained optimization problems. We refer to~\cite{Akcelik:2002a,Akcelik:2006a,Biros:1999a,Biros:2005a,Biros:2005b,Biegler:2003a,Biegler:2007a} for examples.}

\subsection{Formulation}
\label{s:formulation-tumor}

\begin{figure}
\centering
\includegraphics[width=0.9\textwidth]{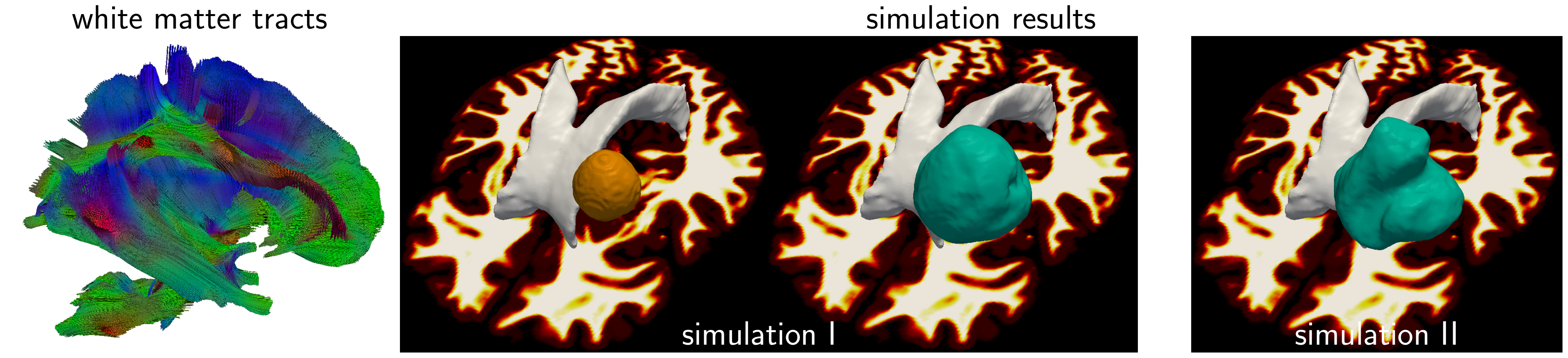}
\caption{Forward simulation. The migration model controls the shape of the tumor by connecting the simulation to virtual brain anatomy (the data is taken from~\cite{Cocosco:1997a}). We display an estimate for white matter fibre tracts of the human brain obtained from diffusion tensor imaging (data taken from~\cite{Mori:2008a}; the color coding encodes the orientation of the fibres). This data allows us to model a preferential migration of cancerous cells along white matter tracts (see, e.g., \cite{Clatz:2005a,Gholami:2016a,Konukoglu:2010a}). We show two simulation results for different model parameters (middle: isotropic migration model; right: anisotropic migration model). The orange surface of the first simulation shows the core of the tumor and the green surface the outer boundary. The images are modified from~\cite{Mang:2014a}.\label{f:forward-simulation}}
\end{figure}

The PDE constraint $\D{C}$ in~\eqref{e:pde-opt}, i.e., the forward operator, is given by
\begin{subequations}
\label{e:tumor:fwd}
\begin{align}
\label{e:tumor:fwd-pde}
\p_t m - \idiv k \igrad m - \rho m(1-m)
& = 0
&& {\rm in}\;\Omega_B\times(0,1],
\\
\label{e:tumor:fwd-initial}
m&=m_I
&& {\rm in}\;\Omega_B\times\{0\},
\end{align}
\end{subequations}

\noindent with Neumann boundary conditions on $\p\Omega_B$. The parameter $\rho > 0$ represents the growth rate; it controls the proliferation of cancerous cells modeled by a logistic function. The diffusion tensor $k \in \fs{M}$, $\fs{M} \defeq \{\tilde{k} :\bar{\Omega}_B \rightarrow \ns{R}^{d,d} \mid\tilde{k}(x) \succ 0 \text{ and } \tilde{k}(x) = \tilde{k}(x)^\T\text{ for all }x\}$, controls the net migration of cancerous cells $m:\tilde{\Omega}_B\times[0,1]\rightarrow[0,1]$. A common assumption in many models is that the rate of migration depends on the tissue composition; typically, it is modeled to be faster in white than in gray matter~\cite{Le:2015a,Le:2017a,Swanson:2000a}. If we denote the probability maps for these tissue types by $\pi_W:\bar{\Omega}\rightarrow[0,1]$ and $\pi_G:\bar{\Omega}\rightarrow[0,1]$, respectively, the associated model for $k$ is given by~\cite{Mang:2014a}
\begin{equation}
\label{e:diffusion-isotropic}
k(x) = (k_W\pi_W(x) + k_G\pi_G(x)) I, \qquad k_W, k_G > 0,
\end{equation}

\noindent with $I \defeq \operatorname{diag}(1,\ldots,1)\in\ns{R}^{d,d}$. We can augment this model by integrating the white matter tract architecture\footnote{The white matter fibre architecture can be estimated from data measured by so-called diffusion tensor imaging, a magnetic resonance imaging technique that measures the anisotropy of diffusion in the human brain. The result of this measurement is a tensor field $\tilde{k}$ that can directly be inserted into~\eqref{e:tumor:fwd}.} (i.e., information about neuronal pathways; see \figref{f:forward-simulation} for an illustration). This has, e.g., been considered in~\cite{Clatz:2005a,Gholami:2016a,Konukoglu:2010a,Mang:2014a,Mang:2012b,Mosayebi:2012a}.

The input to our inversion are \bipa\item estimates for $\pi_W$ and $\pi_G$, \item a tumor free patient image obtained from medical imaging data~\cite{Gooya:2012a,Menze:2015a}, and \item an estimate for the patient's tumor given in terms of a probability map $\pi_T : \bar{\Omega}_B\rightarrow[0,1]$\eipa. Given the forward model~\eqref{e:tumor:fwd} our parameter vector $w$ essentially consists of $\rho$, $k$ (i.e., $k_W$ and $k_G$ in~\eqref{e:diffusion-isotropic}) and $m_I$. We assume that we have experimental data that provides estimates for $k$ and $\rho$; we only invert for $m_I$. We parameterize the initial condition $m_I$ in~\eqref{e:tumor:fwd-initial} using an $n_p$ dimensional space spanned by a Gaussian basis \[\phi_l(x)= (2\pi|\Sigma|^{\half{1}})^{-1}\exp\big(-\half{1}(x - \mu_l)^\T\Sigma^{-1}(x - \mu_l)\big);\] we obtain $m_I(x) \defeq \sum_{l=0}^{n_p} p_l \phi_l(x) = \Phi p$. The regularization model for our problem formulation (see~\eqref{e:pde-opt}) is the $\ell^2$-norm of $p$.

\section{Cardiac Motion Estimation}
\label{s:heart}

Medical imaging can help in the diagnosis of cardiac masses, cardiomyopathy, myocardial infarction, and valvular disease. One example is cine-MRI, which is used routinely for diagnosing a variety of cardiovascular disorders~\cite{Castillo:2003a,Shah:2005a}. Given the cine-MRI data we seek to reconstruct the motion of the cardiac wall, which in turn is used to evaluate clinical indices like ejection fraction and localized abnormalities of the heart function. Segmentation of the ventricles and the myocardium is the first step toward motion reconstruction for quantitative functional analysis of cine-MRI data. Manual segmentation is time consuming~\cite{Axel:2002a} and inter- and intra-observer variability can be high~\cite{Bluemke:2003a}. For this reason automatic segmentation algorithms are preferable. However, even today most of them are not robust enough for a fully automatic segmentation. Here we discuss methods that incorporate biophysical constraints and often result in PDE-constrained optimization algorithms.

One of the main thrusts in recent research has been the 4D motion estimation using biomechanical models. In~\cite{Hu:2003a}, a piecewise linear composite biomechanical model was used to determine active forces and the strains in the heart based on tagged MRI information. In \cite{Papademetris:2001a,Papademetris:2002a} echocardiography and MR images were combined with biomechanical models for cardiac motion estimation. Interactive segmentation was combined with a Bayesian estimation framework that was regularized by an anisotropic, passive, linearly elastic, myocardium model. The authors recognized the importance of neglecting active contraction of the left ventricle. In \cite{Sermesant:2006a,Sermesant:2006b}, the need for realistic simulations and the need for inversion and data assimilation was outlined. In~\cite{Bistoquet:2006a}, an incompressible deformable registration was used for global indicators. In~\cite{Perperidis:2005a}, a regular grid with B-splines was used with a mutual information-based similarity measure;\footnote{Mutual information is a statistical distance measure that originates from information theory. As opposed to the squared $L^2$ distance used in the present work (see \secref{s:diffeomorphic-registration}), which is used for registering images acquired with the same modality, mutual information is used for the registration of images that are acquired with different modalities (e.g., the registration of computed tomography and magnetic resonance images). Mutual information assess the statistical dependence between two random variables (in our case image intensities; see~\cite{Modersitzki:2004a,Modersitzki:2009a,Sotiras:2013a} for more details).} however, the focus has been on inter-subject registration~\cite{Perperidis:2005b} (without consideration of mechanics). In ~\cite{Shen:2005a}, a 4D deformable registration was proposed with a diffusive space-time regularization that worked well for global motion measures. In~\cite{Sundar:2009a}, a biomechanically-constrained motion estimation algorithm that explicitly couples raw image information with a biomechanical model of the heart is described. The main features of that scheme are \bipa\item a patient-specific image-based inversion formulation for the active forces; \item a multigrid-accelerated, octree-based, adaptive finite-element forward solver that incorporates anatomically-accurate fiber information; \item an adjoint/Hessian-based inversion algorithm; and \item a 4D coupled inversion for all images\eipa. The method requires segmentation of the blood pool and myocardium at the initial frame to assign material properties and a deformable registration to map fiber-orientations to the patient. More complex models are considered in~\cite{Delingette:2012a}. The authors use an adjoint method for cardiac-motion estimation and contractility combining catheterized electrophysiology data and cine-MRI sequences. They consider a gradient descent method that uses adjoint equations for the elastic deformation of the heart. A similar approach for cardiac motion estimation is described in~\cite{Tuyisenge:2016a}; the adjoints are derived using an automatic differentiation framework.

\section{Numerics}
\label{s:numerics}

Here, we summarize the numerical strategies for solving the problems in~\secref{s:diffeomorphic-registration} and~\secref{s:dataassim-tumor}, respectively. We discuss the discretization in space and time, the solvers, the computational kernels, and the parallel implementation. Overall, we use a globalized, preconditioned, inexact, reduced space Gauss--Newton--Krylov method for both problems. We will see that the optimality systems of the considered problems yield large, ill-conditioned operators that are challenging to solve in an efficient way. The associated PDE operators will determine the choices we make for the numerics.

\subsection{Optimality Conditions}
\label{s:optimality-systems}

We use the method of Lagrange multipliers~\cite{Lions:1971a} to turn~\eqref{e:pde-opt} into an unconstrained problem. From thereon we can consider two alternative approaches to tackle the optimization problem: \emph{optimize-then-discretize} and \emph{discretize-then-optimize} (see \cite[p.~57ff.]{Gunzburger:2003a,Wilcox:2015a} for a discussion; see also \remref{r:opt-dis-dis-opt}). We use the former strategy; we first compute variations of the Lagrangian functional with respect to the state, adjoint, and control variables, and then discretize them. We will assume that all variables fulfill the necessary regularity requirements to be able to carry out these computations; we refer the interested reader to, e.g.,~\cite{Barbu:2016a,Chen:2011b,Chen:2011a,Lee:2010a,Lee:2011a}, for a rigorous treatment.

\begin{remark}\label{r:opt-dis-dis-opt}
An optimize-then-discretize approach does not guarantee that the discretization of the gradient is consistent with the discretization of the objective functional. Further, it is not guaranteed that the discretization of the forward and adjoint operators are transposes of one another; similarly, it is not guaranteed that discretization of the Hessian is symmetric. These inconsistencies can have a negative effect on the convergence of the solver. For a discretize-then-optimize approach one can (by construction) guarantee that these operators are consistent, transposes, and symmetric. However, one may, e.g., obtain a different numerical accuracy for the forward and adjoint operators (see, e.g.,~\cite{Hager:2000a,Dontchev:2000a} for examples). For the optimize-then-discretize approach we can choose arbitrarily accurate solvers for the forward and adjoint problems. We refer to~\cite{Gunzburger:2003a,Borzi:2012a} for more remarks on the discretization of optimization and control problems. For the PDE-constrained optimization problems considered in this study, we have experimentally observed that the discretization errors are below the tolerances relevant in practical application (e.g., a relative change of the gradient of $\num{1E-2}$). If we are interested in solving the inverse problems more accurately (e.g., up to $\num{1E-6}$) our solver may fail to converge. If more accurate solutions are required, we will have to consider PDE solvers that preserve the properties of the continuous operators. For the diffeomorphic registration problem, we present a discretize-then-optimize implementation in~\cite{Mang:2017a}; works of other authors that consider a discretize-then-optimize strategy can be found in~\cite{Ashburner:2011a,Benzi:2011a,Herzog:2018a}.
\end{remark}

\subsubsection{Diffeomorphic Registration}
\label{s:optimality-systems:reg}

We assume that $v$ is stationary and divergence-free. We introduce the Lagrange multipliers $\lambda:\bar{\Omega}\times[0,1]\rightarrow\ns{R}$, $\lambda_0:\bar{\Omega}\rightarrow\ns{R}$, $\mu:\bar{\Omega}\rightarrow\ns{R}$, to obtain the Lagrangian functional\footnote{We neglect the periodic boundary conditions for simplicity.}
\begin{align}
\label{e:lagrangian-reg}
\F{L}_R[m,w,v,\lambda,\lambda_0,\mu]
&\defeq
\half{1}\|m_1-m_R\|_{L^2(\Omega)}^2 + \half{\beta}\|v\|_{\fs{V}}^2
+ \langle \idiv v, \mu\rangle_{L^2(\Omega)} \\\nonumber
& + \int_0^1\langle \p_t m + \igrad m \cdot v,\lambda \rangle_{L^2(\Omega)}\d{t}
+ \langle m_0 - m_T,\lambda_0\rangle_{L^2(\Omega)}
\end{align}

\noindent We obtain the optimality system by taking variations with respect to the adjoint, state, and control variables, and applying integration by parts. We use a reduced space formulation. That is, we assume that the state and adjoint equation are fulfilled exactly (see below) and only iterate on the reduced space of the velocity $v$. Accordingly, we require vanishing variations of the Lagrangian with respect to the control variable $v$, i.e., \[g_v(v) \overset{!}{=} 0\quad\text{in}\; \Omega.\] This is the strong form of the control or decision equation (i.e., the reduced gradient) of our problem; $g_v(v)$ is an integro-differential operator of the form
\begin{align}
\label{e:control}
g_v(v) \defeq \beta\D{A}[v]
+ \D{K}[\int_0^1\lambda\igrad m\d{t}\;].
\end{align}

\noindent The operator $\D{A}$ originates from the regularization operator in~\eqref{e:reg:regularization}; depending on its order we arrive at an elliptic~\cite{Mang:2016a}, biharmonic~\cite{Mang:2015a} or triharmonic~\cite{Chen:2011a} integro-differential control equation. The pseudo-differential operator $\D{K}$ is a projection onto the space of divergence-free velocities. It originates from the elimination of the constraint $\idiv v = 0$ and the adjoint variable $\mu$ from our optimality system (see \cite{Mang:2015a,Mang:2016a}). The operator is given by $\D{K} \defeq I - \igrad\ilap^{-1} \idiv$ (\emph{Leray projection}). If we neglect the incompressibility constraint in our formulation $\D{K}$ becomes an identity operator. If we introduce a mass-source map $q$ into our formulation this operator becomes more involved~\cite{Mang:2016a}. Evaluating the gradient~\eqref{e:control} for a given trial $v$ requires us to find $m$ and $\lambda$. We can compute $m$ by solving the forward problem~\eqref{e:reg:fwd} \emph{forward in time}. Once we have found the \emph{final state} $m_1(x)\defeq m(x,t=1)$ we can compute $\lambda$ by solving the adjoint equation
\begin{subequations}
\label{e:reg:adjoint}
\begin{align}
-\p_t \lambda - \idiv v\lambda & = 0
&&\text{in}\;\; \Omega \times [0,1),
\label{e:reg:adjoint-pde}
\\
\lambda &= m_R - m
&&{\rm in}\;\; \Omega \times\{1\},
\label{e:reg:adjoint-final}
\end{align}
\end{subequations}

\noindent with periodic boundary conditions on $\p\Omega$ \emph{backward in time}; \eqref{e:reg:adjoint} is a final value problem and models the transport of the residual difference between the final state and the reference image backward in time.\footnote{Instead of solving an advection or continuity equation as done in the Eulerian formulation we present here, it has been suggested to solve for the state and adjoint variables using the diffeomorphism $y$ (which involves interpolation operations instead of the solution of a transport equation)~\cite{Vialard:2012a}.}

We can use the expression for the reduced gradient in~\eqref{e:control} within a gradient descent scheme for numerical optimization (see, e.g., \cite{Mang:2015a}). We use a Newton--Krylov method instead. We provide additional details in~\secref{s:gn-krylov-method} and~\secref{s:newton-step}.

\subsection{Data Assimilation}
\label{s:optimality-systems:tumor}

The Lagrangian for our problem is given by\footnote{We neglect the Neumann boundary conditions for simplicity.}
\begin{align}
\label{e:lagrangian-tumor}
\F{L}_T[m,p,\lambda,\lambda_0]
&\defeq
\half{1}\|Qm-\pi_T\|_{L^2(\Omega_B)}^2
+ \half{\beta} \|p\|^2_2
+ \langle m_0 -\Phi p,\lambda_0\rangle_{L^2(\Omega_B)}
\\\nonumber
&+
\int_0^1\langle
\p_t m - \idiv k \igrad m - \rho m(1-m),\lambda
\rangle_{L^2(\Omega_B)}\d{t}\nonumber
\end{align}

\noindent with Langrange multiplier $\lambda:\bar{\Omega}_B \times [0,1] \rightarrow \ns{R}$, and $\lambda_0(x) = \lambda(x,t=0)$, $\lambda_0:\bar{\Omega}_B\rightarrow\ns{R}$, and $m_0(x) \defeq m(x,t=0)$. The unknowns of our problem are the diffusion coefficients $k_W$ and $k_G$, the growth rate $\rho$, and the initial condition $m_I$ (and its location). We can derive the adjoint equations for either or all of these variables. For simplicity of presentation, we will derive the optimality conditions for the parameterization $p$ of $m_I$ in~\eqref{e:tumor:fwd-initial}. A more complete picture can be found in~\cite{Gholami:2016a}. The strong form of the control or decision equation (i.e., the reduced gradient of our problem) is given by the vanishing first variation of the Lagrangian functional $g_p(p) = 0$ in $\ns{R}^{n_p}$, where
\begin{align}
\label{e:tumor:gradient}
g_p(p) \defeq \beta p - {\Phi}^\T\lambda_0.
\end{align}

\noindent To obtain $\lambda_0$, we have to solve the adjoint equation
\begin{subequations}
\label{e:tumor:adj}
\begin{align}
- \partial_t \lambda
- \idiv k \igrad\lambda
- \rho(1- 2m)\lambda
& = 0
&&{\rm in}\; \Omega_B \times [0,1),
\label{e:tumor:adj-pde}
\\
 \lambda
& = \D{Q}^\ast(\D{Q}m - \pi_T)
&&{\rm in}\; \Omega_B \times \{1\},
\label{e:tumor:adj-final}
\end{align}
\end{subequations}

\noindent with Neumann boundary conditions on $\p\Omega_B$, backward in time. Notice that the final condition in~\eqref{e:tumor:adj-final} depends on the final state $m$ at $t=1$, i.e., requires the solution of the forward problem~\eqref{e:tumor:fwd}.

\subsection{Discretization}
\label{s:discretization}

We use the trapezoidal rule for numerical quadrature with respect to space and with respect to time. We assume that the functions in our formulation (including images) are periodic and continuously differentiable.\footnote{The data assimilation problem in~\secref{s:dataassim-tumor} requires Neumann boundary conditions on $\p\Omega_B$, with $\Omega_B\subset\Omega$. We use a penalty approach to approximate these boundary conditions. We apply periodic boundary conditions on $\p\Omega$ and set the reaction and diffusion coefficients in~\eqref{e:tumor:fwd-pde} to sufficiently small penalty parameters $k^\epsilon\rightarrow 0$ and $\rho^\epsilon\rightarrow0$ outside of $\Omega_B$; see, e.g.,~\cite{Gholami:2016a,Hogea:2008b,Mang:2014a,Mang:2012a}.}  We apply appropriate filtering operations and periodically extend or mollify the discrete data to meet these requirements. We use regular grids to discretize $\Omega\times[0,1]$, $\Omega \defeq [0,2\pi)^d\subset\ns{R}^d$, $d\in\{2,3\}$. The spatial grid consists of $n = \textstyle\prod_{i=1}^d n_i$, $n_i\in\ns{N}$, grid points $\di{x}_j\defeq(x_{1,j},\ldots,x_{d,j})^\T\in\ns{R}^d$, $x_{i,j}\defeq2\pi j_i/n_i$, $0 \leq j_i \leq n_i-1$, $i=1,\ldots,d$. We consider a nodal discretization in time, which results in $n_t+1$ discretization points. We use a spectral projection scheme for all spatial operations. That is, we approximate spatial functions as $\textstyle m(x)=\sum_k \hat{m}_k\exp(-k\cdot x)$, where $k=(k_1,\ldots,k_d)^\T\in\ns{N}^d$, $-n_i/2+1 \leq k_i \leq n_i/2$, $i=1,\ldots,d$, is the grid index and $\hat{m}_k$ are the spectral coefficients of $m$. The mappings between the spatial and spectral coefficients are done using FFTs. We approximate spatial derivative operators by applying the appropriate weights in the spectral domain. This scheme allows us to efficiently, stably, and accurately apply differential operators and their inverses.\footnote{Note that the $\inf$-$\sup$ condition for pressure spaces arising in finite element discretizations of Stokes problems~\cite[p.~200ff.]{Brezzi:1991a} is not an issue with our scheme (incompressibility constraint in diffeomorphic registration problem).}

\subsection{Numerical Time Integration}
\label{s:time-integration}

The iterative solution of the optimality conditions of our problems requires a repeated solution of PDE operators of mixed type. Here we summarize the time integration schemes for the parabolic or hyperbolic PDEs that appear in our optimality systems. In general we opt for unconditionally stable schemes to reduce the number of timesteps. This will not only reduce the computational work load but also reduce the memory requirements; evaluating the Hessian operators that appear in our scheme requires us to store at least one space-time field (the state variable).\footnote{A possible remedy is to employ check-pointing or domain decomposition strategies~\cite{Akcelik:2002a,Griewank:1992a,Heinkenschloss:2005a}.}

\subsubsection{Parabolic PDEs}
\label{s:time-integration:parabolic}

We use an unconditionally stable, second-order Strang-splitting method to solve the parabolic equations (see also \cite{Gholami:2016a,Hogea:2008a}). We explain this for the forward problem~\eqref{e:tumor:fwd} in \algref{a:strang-splitting}. Let $\di{m}^j\in\ns{R}^n$ denote the tumor discretized distribution at time $t^j = jh_t$, $h_t = 1/n_t$. We apply an implicit Crank--Nicolson method for diffusion (lines~\ref{lst:line:diff1} and~\ref{lst:line:diff2} in~\algref{a:strang-splitting}) and solve the reaction part analytically (line \ref{lst:line:prol} in \algref{a:strang-splitting}). A similar splitting strategy is used in~\cite{HawkinsDaarud:2013a}. In other related work, conditionally~\cite{Colin:2014a,Weis:2017a} and unconditionally stable~\cite{Mang:2014a,Mang:2012a} explicit schemes were considered.

\begin{algorithm}
\caption{Strang-splitting for the solution of the forward problem in \eqref{e:tumor:fwd}. We denote the discretized diffusion operator by $\di{L}=[\idiv]^h \di{k}\igrad^h$.}
\label{a:strang-splitting}
\begin{algorithmic}[1]
\STATE{input: $h_t \leftarrow 1/n_t$, $p\in\ns{R}^{n_p}$, $\di{\Phi} \in \ns{R}^{n,n_p}$}
\STATE{$\di{m}^0\leftarrow\di{\Phi}p$; $\di{I} \leftarrow\operatorname{diag}(1,\ldots,1)\in\ns{R}^{n,n}$}
\FOR{$j=0,\ldots,n_t$}
\STATE{$\di{m}^\dagger \leftarrow$ solve $(\di{I} - \frac{h_t}{4}\di{L})\di{m}^\dagger = (\di{I} + \frac{h_t}{4}\di{L})\di{m}^j$ in $(t^j,t^j+\half{h_t}]$\label{lst:line:diff1}}
\STATE{$\di{m}^\ast\leftarrow$ solve $\p_t \di{m}= \rho\di{m}(1-\di{m})$ in $(t^j,t^{j+1}]$ analytically with initial condition $\di{m}^\dagger$\label{lst:line:prol}}
\STATE{$\di{m}^{j+1} \leftarrow$ solve $(\di{I} - \frac{h_t}{4}\di{L})\di{m}^{j+1} = (\di{I} + \frac{h_t}{4}\di{L})\di{m}^\ast$ in $(t^j+\half{h_t},t^{j+1}]$\label{lst:line:diff2}}
\ENDFOR
\end{algorithmic}
\end{algorithm}

We use a preconditioned conjugate gradient (PCG) method with a fixed tolerance of $\num{1E-6}$ to solve the diffusion steps. We use a preconditioner $\di{P}$ based on a constant coefficient approximation of $\di{L}$ given by $\di{P} = \di{I} - \frac{h_t}{4}\tilde{k}\ilap^{\!h}$, where $\tilde{k} > 0$ is the average diffusion coefficient. The inversion and construction of $\di{P}$ has vanishing computational cost due to our spectral scheme; it only requires a pointwise multiplication in the Fourier domain and a single forward and backward FFT.

\subsubsection{Hyperbolic PDEs}
\label{s:time-integration:hyperbolic}

We employ a semi-Lagrangian scheme~\cite{Falcone:1998a} to solve the hyperbolic transport equations of the form $\p_t m + \igrad m \cdot v = f(m,v)$. In the context of diffeomorphic image registration, semi-Lagrangian schemes have, e.g., been used in~\cite{Beg:2005a,Chen:2011a,Mang:2017c,Mang:2016c}. They are unconditionally stable. This allows us to keep the number of timesteps $n_t$ small and by that significantly reduces the memory requirements. We illustrate this scheme for the forward problem (homogeneous case) in~\algref{a:semi-lagrangian}; additional details can be found in~\cite{Mang:2017c}. The algorithm consists of two steps: In a first step we have to compute the characteristic $\di{y}$ by solving $\d_t^h \di{y} = \di{v}(\di{y})$ in $[t^j,t^{j+1})$ with final condition $\di{y} = \di{x}$ at $t^{j+1}$ backward in time. The second step is to solve an ODE for the transported quantity along $\di{y}$; i.e., we have to solve $\d_t^h \di{m}(\di{y}) = \di{f}(\di{m},\di{v}; \di{y})$. Both of these steps require the evaluation of scalar or vector fields along the characteristic $\di{y}$, i.e., at locations that do not coincide with grid points $\di{x}$; this involves interpolation (see lines \ref{lst:interpol1} and \ref{lst:interpol2} in \algref{a:semi-lagrangian} for an example). We use an explicit, second-order accurate Runge--Kutta scheme to integrate these ODEs in time, and a cubic interpolation model in space.

\begin{algorithm}
\caption{Semi-Lagrangian scheme for the solution of the forward problem in~\eqref{e:reg:fwd}. Note, that we have $\d_t^h \di{m}(\di{y}) = \di{0}$ for the forward problem.}
\label{a:semi-lagrangian}
\begin{algorithmic}[1]
\STATE{input: $h_t \leftarrow 1/n_t$, $\di{m}_T \in \ns{R}^n$, $\di{x}\in\ns{R}^{dn}$, $\di{v}\in\ns{R}^{dn}$}
\STATE{$\di{\tilde{y}}$ $\leftarrow$ $\di{x} - h_t \di{v}$}
\STATE{$\di{y}$ $\leftarrow$ $\di{x} - \half{h_t} ({\tt interpolate}(\di{v},\di{\tilde{y}}) + \di{v})$ \label{lst:interpol1}}
\STATE{$\di{m}^0\leftarrow\di{m}_T$}
\FOR{$j=0,\ldots,n_t$}
\STATE{$\di{m}^{j+1}$ $\leftarrow$ ${\tt interpolate}(\di{m}^j,\di{y})$ \label{lst:interpol2}}
\ENDFOR
\end{algorithmic}
\end{algorithm}

\subsection{Numerical Optimization}

We have derived the expressions for the optimality conditions in~\secref{s:optimality-systems}. Here, we revisit the generic formulation used in~\eqref{e:pde-opt}. The specific operators for the two problems can be found in \secref{s:optimality-systems} and \secref{s:newton-step}.

In general, we require that the gradient with respect to the control variable $\di{w}$ vanishes; i.e., $\di{g}(\di{w}^{\star}) = 0$ for an admissible solution $\di{w}^{\star}$ to~\eqref{e:pde-opt}. We use a globalized, inexact~\cite{Dembo:1983a,Eisenstat:1996a}, preconditioned, matrix-free, reduced space\footnote{By \emph{reduced space} we mean that we iterate only on the reduced space of the control variable of our problem as opposed to so called \emph{full-space} or \emph{all-at-once} methods (see, e.g.,~\cite{Benzi:2009a,Biros:2005a,Biros:2005b,Haber:2001a,Herzog:2010a} for more details).} Gauss--Newton--Krylov method for numerical optimization. The Newton step is (in general format) given by
\begin{equation}
\label{e:reduced-space-kkt-system}
\di{w}_{k+1} = \di{w}_k + \alpha_k\tilde{\di{w}}_k, \quad k = 1,2,\ldots
\end{equation}

\noindent where $\di{w}_k, \tilde{\di{w}}_k\in\ns{R}^{\tilde{n}}$.\footnote{The order $\tilde{n}$ of the optimality system depends on the problem. For the diffeomorphic registration case the control variable $\di{w}$ is given by the velocity field $\di{v}\in\ns{R}^{dn}$, i.e., $\tilde{n}\equiv dn$; for the tumor case $\di{w}$ is given by $p\in\ns{R}^{n_p}$, i.e., $\tilde{n}\equiv n_p$.} We compute the search direction $\tilde{\di{w}}_k$ by solving the reduced space Karush--Kuhn--Tucker (KKT) system $\di{H}_k \tilde{\di{w}}_k = -\di{g}_k$ iteratively using a PCG method. Here, $\di{H}_k \in \ns{R}^{\tilde{n},\tilde{n}}$, $\tilde{n}\in\ns{N}$, $\di{H}_k \succ 0$, is a discrete representation of the Gauss--Newton approximation (see \remref{r:gn-approx}) to the reduced Hessian at Newton (or outer) iteration index $k$. We globalize our scheme based on a backtracking line search subject to the Armijo condition (we use default parameters; see \cite[algorithm~3.1,\,p.~37]{Nocedal:2006a}).\footnote{Our implementation also features a trust region method.}

Since we use a matrix-free Krylov-subspace method to iteratively solve for the search direction $\tilde{\di{w}}_k$, we only require an expression for the application of the Hessian to a vector (\emph{Hessian matvec}). In Newton-type methods, most work is spent on the (iterative) solution of the KKT system; it is the computational bottleneck of our solver. For both of our problems, each Hessian matvec involves the solution of two PDEs: the incremental state and adjoint equations; this is costly. Therefore, it is essential to design an effective preconditioner to reduce the number of Hessian matvecs we have to perform. We summarize the overall algorithm in~\algref{a:outer-iteration} (Newton step) and the steps necessary to solve~\eqref{e:reduced-space-kkt-system} in~\algref{a:inner-iteration} (see \secref{s:gn-krylov-method}). We provide expressions for the Hessian matvec for the two individual problems in~\secref{s:newton-step}, respectively.

The preconditioner for the reduced space Hessian of the tumor problem uses, likewise to the preconditioner for the forward problem described \secref{s:time-integration:parabolic}, a constant coefficient approximation for the diffusion operators (see~\cite{Gholami:2016a}). For the registration problem we use the inverse of the regularization operator $\D{A}$ as a preconditioner (a common choice in PDE-constrained optimization~\cite{Alexanderian:2016a,Mang:2015a,Mang:2016c,Mang:2017a}). The preconditioned Hessian is a compact perturbation of identity; for smooth data and sufficient regularity of $v$ (i.e., we can resolve the problem), we observe a rate of convergence that is independent of the mesh size. In recent work, we have presented a 2-level (multigrid inspired) preconditioner. The convergence for both preconditioners is not $\beta$-independent. Details for the 2D case can be found in~\cite{Mang:2017c}. The results presented in this study are in 3D and computed using $\D{A}^{-1}$ as a preconditioner.

\begin{remark}\label{r:gn-approx}
We use a Gauss--Newton approximation to the Hessian to ensure that the operator is positive definite far away from the optimum. The Gauss--Newton approximation is obtained by dropping all expressions in the evaluation of the Hessian that involve the adjoint variable $\digr{\lambda}$. (We provide a precise definition of the individual operators in~\secref{s:newton-step}.) We expect that the rate of convergence of our scheme drops from quadratic to superlinear using this approximation. We recover local quadratic convergence as $\digr{\lambda}\rightarrow \di{0}$ (i.e., the mismatch goes to zero, and we have solved the problem). We observe this for synthetic problems. For problems involving real data, the data will be perturbed by noise. In fact, we may not only have to deal with noise perturbations, but also with intensity drifts (which we can correct for in a preprocessing step). Based on these perturbations, we may not recover local quadratic convergence.
\end{remark}

\subsection{Implementation Details}
\label{s:implementation-details}

Generally speaking, (interesting) images are functions of bounded variation. Our scheme cannot handle this type of discontinuities in the data. We assume that all fields that appear in our formulations are smooth, compactly supported functions. We ensure this numerically by presmoothing the data, and applying appropriate extensions to the data, if necessary. For the registration problem we estimate the regularization parameter based on a parameter continuation scheme~\cite{Haber:2000a,Haber:2006a,Mang:2015a}. For the tumor problem we use a similar approach---the L-curve method~\cite{Hansen:1992a} (see~\cite{Gholami:2016a} for an example). We terminate the inversion if the relative change of the gradient of our problem reaches a user defined tolerance. We typically choose values between $\num{1e-1}$ and $\num{1E-2}$ for real world applications.

\begin{remark}
From a numerical optimization point of view, we want to drive the gradient to zero. In diffeomorphic image registration problems on real data (mainly multi-subject neuroimaging applications) we have observed that, while the solution may still change, the mismatch seems to be quite stable after reducing the gradient by roughly one order of magnitude. Reducing the gradient further does, in our experience, not dramatically alter the mismatch (which is what practitioners mostly care about).
\end{remark}

We have developed and presented prototype implementations of our nonlinear solver in Matlab. These solvers are described in detail in~\cite{Mang:2015a,Mang:2017c,Gholami:2016a}. We summarize the Newton--Krylov scheme in \secref{s:gn-krylov-method}. Our parallel code~\cite{Mang:2016a,Gholami:2017a} is implemented in {\tt C++}, and uses the message passing interface (MPI) for parallelism. We use a 2D pencil decomposition for 3D FFTs~\cite{Grama:2003a,Czechowski:2012a} to distribute the data among processors. We use PETSc for linear algebra operations~\cite{Balay:2016a}, PETSc's TAO package for the nonlinear optimization~\cite{Munson:2017a}, and AccFFT for Fourier transforms~\cite{Gholami:2016b} (a library for parallel FFTs developed by our group). Our code implements the operations for evaluating the objective functional, the gradient, and the Hessian matvec, and the interfaces to set tolerances and control the stopping criteria. The main computational kernels of our scheme are the FFT used for spatial differentiation and the interpolations in the semi-Lagrangian scheme. These kernels, and their parallel implementation, are described in detail in~\cite{Gholami:2016b,Mang:2016c,Gholami:2017a}.

\section{Numerical Experiments}
\label{s:numerical-experiments}

We showcase some numerical experiments to illustrate the behavior of our solvers next.

\subsection{Diffeomorphic Image Registration}
\label{s:experiments:reg}

A study for the performance of our Gauss--Newton--Krylov scheme and a comparison to (Sobolev) gradient descent type approaches can be found in~\cite{Mang:2015a}. We could demonstrate that our Gauss--Newton--Krylov scheme yields a significantly better rate of convergence. A study on registration quality as a function of different regularization operators can be found in~\cite{Mang:2016a}. In this work, we showed that we can generate high-fidelity registration whilst maintaining well-behaved $\det\igrad y$, by constraining local volume change. A study of the performance of our semi-Lagrangian scheme in combination with a two-level preconditioner be found in~\cite{Mang:2017c}. Here, we could achieve a speedup of up to 20x compared to our original solver~\cite{Mang:2015a}. All this work is limited to the 2D case and implemented in \emph{Matlab}. In this section, we present results for the implementation of a memory-distributed solver for the 3D case. For a study on the scalability of this solver on supercomputing platforms we refer to~\cite{Mang:2016c,Gholami:2017a}. In \cite{Gholami:2017a} we achieved a 8x speedup compared to our original implementation~\cite{Mang:2016c}. With this work we pave the way for solving diffeomorphic image registration problems on real clinical data sets for clinically relevant problem sizes in real-time---with a time to solution of under two seconds for 50 million unknowns (a grid size of $256\times300\times256$) using 512 cores on TACC's Lonestar 5 system (2-socket Intel Xeon E5-2690 v3 (Haswell) with 12 cores/socket, \SI{64}{\giga\byte} memory per node). We were also able to solve an inverse problem of unprecedented scale with $\sim$200 million unknowns (a problem size, that is 64x bigger than in~\cite{Mang:2016c}) using 8192 cores on HLRS's Hazel Hen system (2-socket Intel Xeon E5-2690 v3 (Haswell) with 12 cores/socket, \SI{128}{\giga\byte} memory per node; we summarize these results in \tabref{s:scaling-reg}).

\begin{table}
\centering
\caption{Scaling results for diffeomorphic image registration. Top block (left): Strong scaling on TACC's Lonestar 5 system for 3D brain imaging data. Bottom block (left): Weak scaling results on HLRS's HazelHen system. We report (from left to right) the number of discretization points in space, the number of cores (MPI tasks), the runtime (in seconds for the entire inversion), and the strong and weak scaling efficiency (in \%). On the right we plot the strong scaling performance. We also show the ideal performance and the time spent on computing FFTs and evaluating the interpolation operator. Results originally reported in~\cite{Gholami:2017a}.\label{s:scaling-reg}}
\begin{minipage}[b]{0.49\textwidth}
\centering
\begin{scriptsize}
\begin{tabular}[b]{rrrr}\toprule
 problem size            & cores & runtime       & efficiency  \\\midrule
 $256\times300\times256$ & 2     & \num{2.358803e+02} & 100.00 \\
                         & 8     & \num{6.730872e+01} &  87.61 \\
                         & 32    & \num{1.811838e+01} &  81.37 \\
                         & 128   & \num{4.633811e+00} &  79.54 \\
                         & 512   & \num{1.523017e+00} &  60.50 \\
\midrule
1024\textsuperscript{3}  &   128 & \num{1.969103e+02} & 100.0  \\
2048\textsuperscript{3}  &  1024 & \num{2.103693e+02} &  93.6  \\
4096\textsuperscript{3}  &  8192 & \num{2.375269e+02} &  82.9  \\
\bottomrule
\end{tabular}
\end{scriptsize}
\end{minipage}
\begin{minipage}[t]{0.4\textwidth}
\centering
\includegraphics[height=3.0cm]{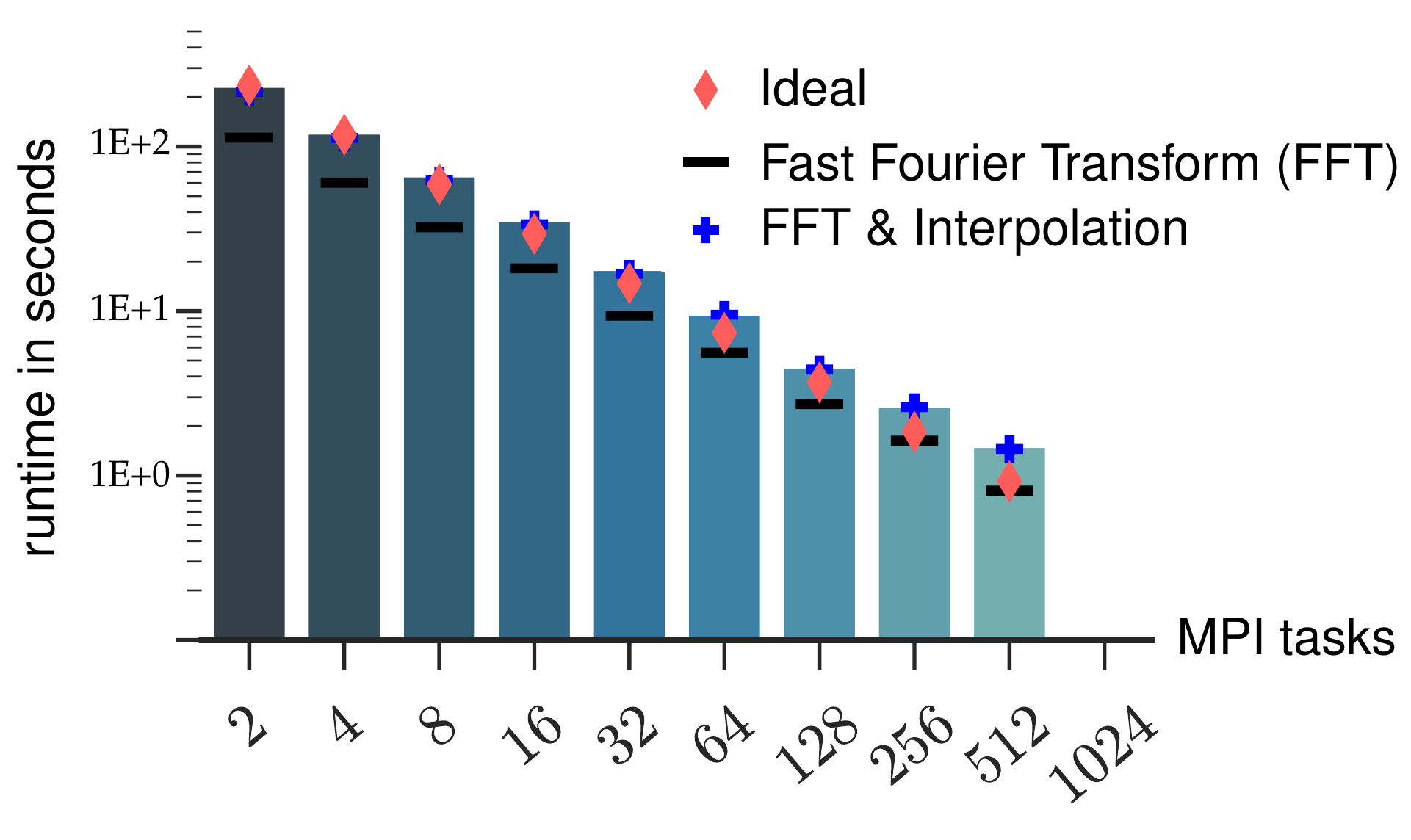}
\end{minipage}
\end{table}

In what follows, we will showcase some new results for our algorithm for diffeomorphic registration~\cite{Mang:2015a,Mang:2016a,Mang:2017c,Mang:2016c}. We will start by examining the capabilities of our solver of modeling large diffeomorphic deformations in 3D. We will then consider synthetic (smooth) and real brain imaging data. We will report results that showcase the performance of our solver on moderate computing platforms.

\subsubsection{Default Parameters}

We normalize the images to an intensity range of $[0,1]$ prior to registration. If we consider data with sharp edges or real images, we apply a Gaussian smoothing operator with standard deviation of one voxel per spatial direction. We perform the registration in double precision. The number of time steps for the PDE solver is fixed to $n_t = 4$. We invert for a stationary velocity field $v$. We use a PCG method to solve for the search direction with a superlinear forcing sequence. We use the inverse of the regularization operator as a preconditioner. The results presented below are obtained on a single node of CACDS's Opuntia cluster (2-socket Intel Xeon E5-2680v2 at \SI{2.8}{\giga\hertz} with 10 cores/socket and \SI{64}{\giga\byte} memory). We do not perform any grid, parameter, or scale continuation. All results reported in this study are computed on the original resolution of the images.

\subsubsection{Sphere-To-Bowl}

We showcase results for a benchmark problem for large deformation diffeomorphic image registration algorithms.

\begin{figure}
\centering
\includegraphics[width=0.6\columnwidth]{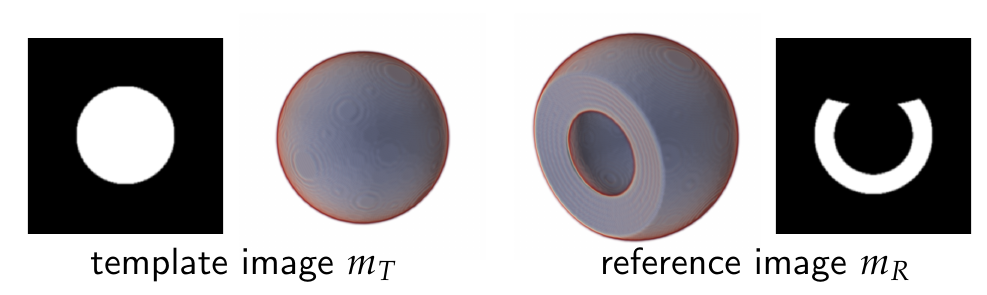}
\caption{Synthetic registration problem as originally proposed in~\cite{Christensen:1996a}. We show two and three-dimensional (volume rendering) illustrations of the template image $m_T$ (left; sphere) and the reference image $m_R$ (right; bowl). The two-dimensional illustration is a plane cut through the center of the three-dimensional objects (sagittal cut; see \figref{f:results-sphere-to-bowl}).}
\label{f:registration-problem:sphere-to-bowl}
\end{figure}

\ipoint{Setup}: We consider the synthetic registration problem proposed in~\cite{Christensen:1996a}---a standard benchmark problem for diffeomorphic image registration. We illustrate the template and reference image in~\figref{f:registration-problem:sphere-to-bowl}. The task is to register a sphere to a bowl. The images are discretized using a grid of size $256^3$ (\num{50331648} unknowns). We use an $H^2$-regularization model and invert for a stationary velocity field $v$.

\begin{figure}
\includegraphics[width=\columnwidth]{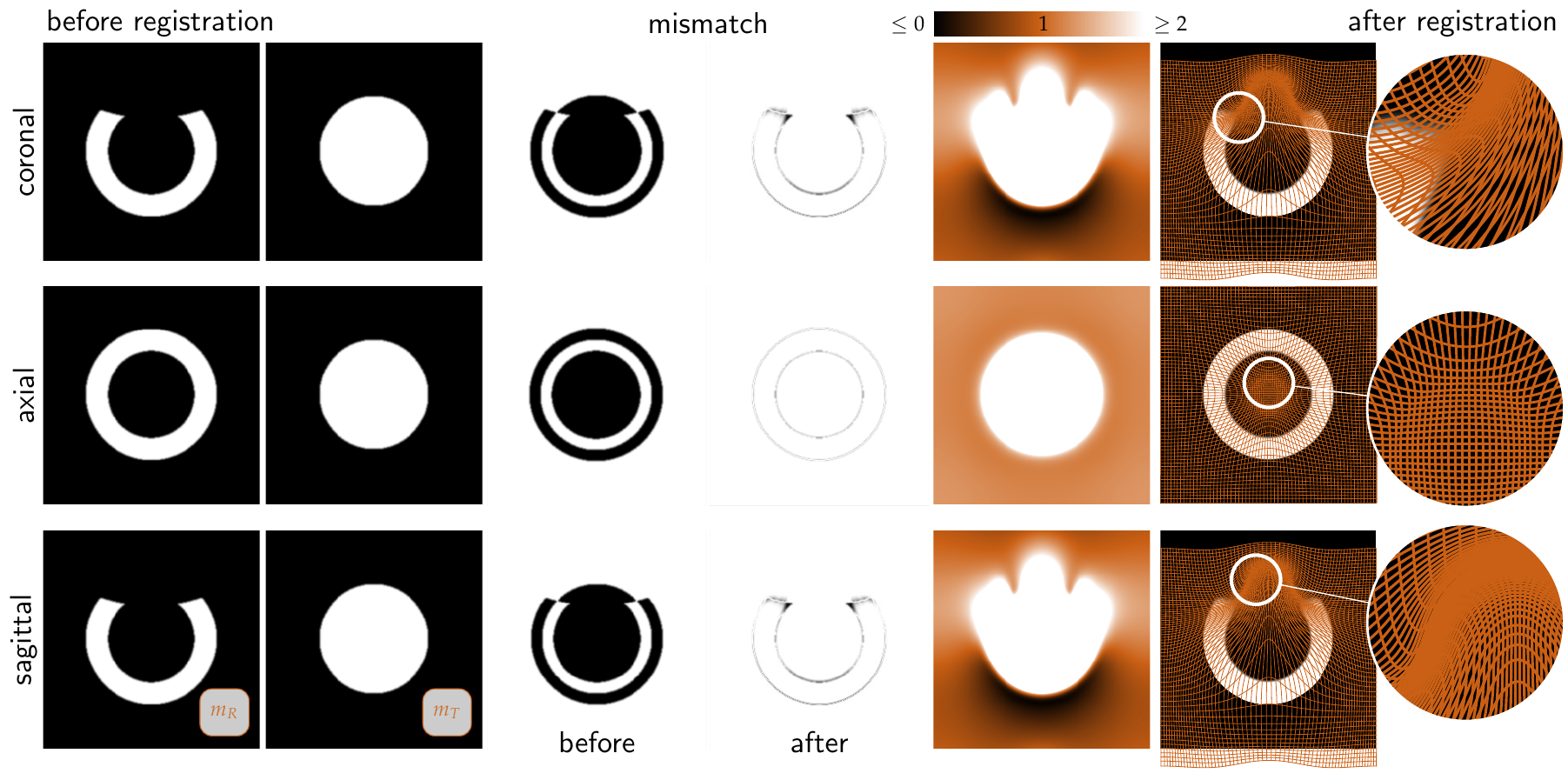}
\caption{Results for the registration problem illustrated in \figref{f:registration-problem:sphere-to-bowl}. We show (from left to right) the reference image $m_R$, the template image $m_T$, the residual differences before registration, the residual differences after registration, a pointwise map of the determinant of the deformation gradient, and the deformed template image with an illustration of the deformed grid in overlay. Each row corresponds to a different view of the 3D volume (from top to bottom: coronal, axial, and sagittal). The map is diffeomorphic as judged from the visualization of the deformed grid as well as the values for the determinant of the deformation gradient (in $[\num{5.302593401e-02},\num{2.732119508e+01}]$; colormap illustrated on the top). Notice that the deformation map illustrated on the left is the pullback map (depicts, where points originate from), whereas the map for the determinant of the deformation gradient corresponds to the inverse of this deformation map.}
\label{f:results-sphere-to-bowl}
\end{figure}

\begin{figure}
\includegraphics[width=\columnwidth]{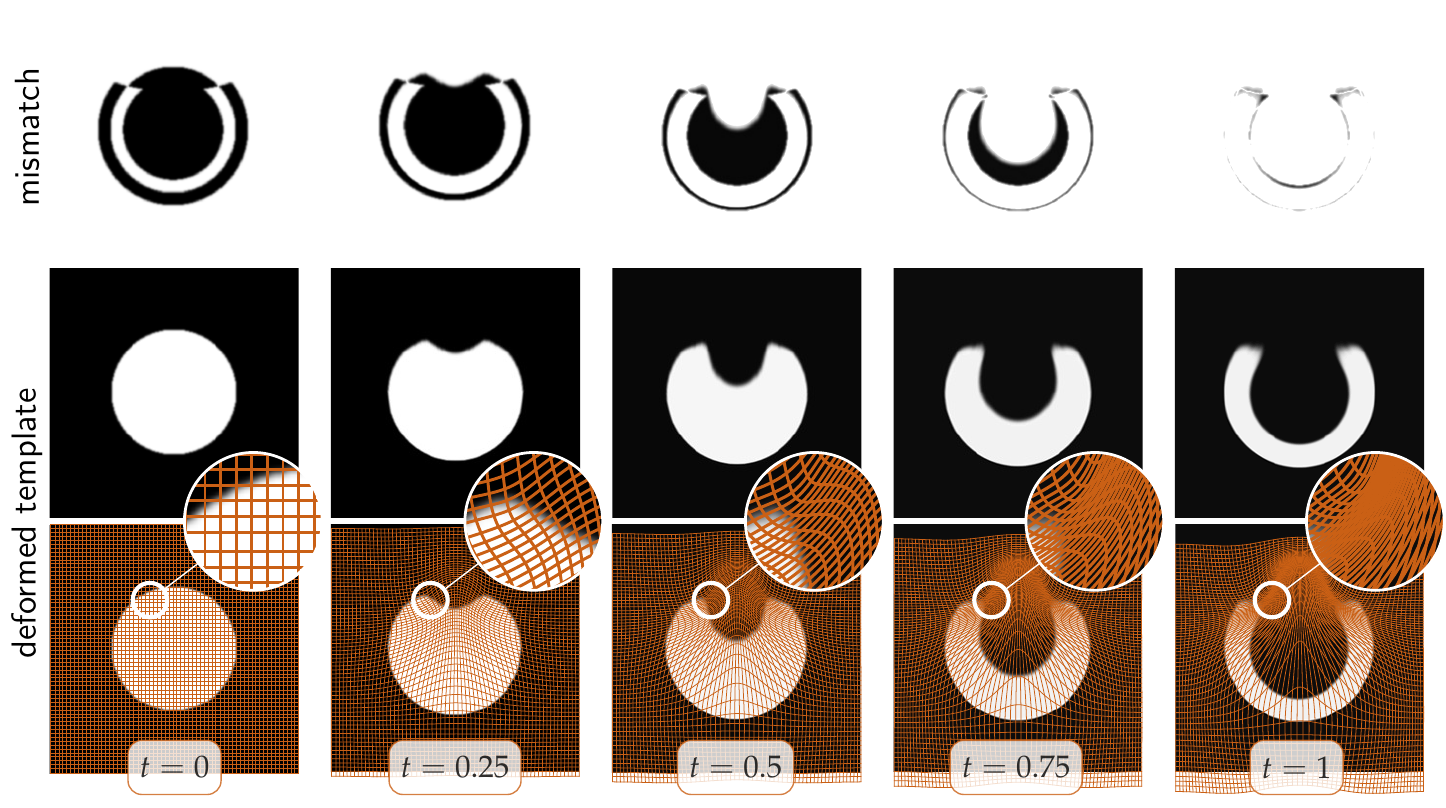}
\caption{Results for the registration problem illustrated in \figref{f:registration-problem:sphere-to-bowl} (sagittal view only). We show (from left to right) the solution of the forward problem for the estimated velocity field $v$ at different time points. The bottom row shows the deformed template image as it evolves in time (bottom row: grid to illustrate the deformation pattern in overlay). The top row shows the mismatch between the deformed template image and the reference image (see \figref{f:registration-problem:sphere-to-bowl} left column for an illustration of the reference image). White areas indicate a vanishing mismatch and black areas a large mismatch.}
\label{f:results-sphere-to-bowl-forwardsolve}
\end{figure}

\ipoint{Results}: We illustrate the results in~\figref{f:results-sphere-to-bowl}. We show three orthogonal planes through the center of the 3D volumes in~\figref{f:registration-problem:sphere-to-bowl} (top to bottom: coronal, axial and sagittal views). Each column in~\figref{f:results-sphere-to-bowl} shows (from left to right) \bipa\item the reference image, \item the template image, \item the mismatch between the reference image and the template image before registration, \item the mismatch between the template and the reference image after registration, \item the determinant of the deformation gradient, i.e., the determinant of the Jacobian of the deformation map $y$, and \item an illustration of the deformation pattern\eipa. The mismatch is large in black areas (corresponds to 100\%) and small in white areas (corresponds to 0\%). The color map for the determinant of the deformation gradient is illustrated in \figref{f:results-sphere-to-bowl}. Black corresponds to a $\det \igrad y = 0$, orange to $\det \igrad y \approx 1$ and white to $\det \igrad y \geq 2$. In \figref{f:results-sphere-to-bowl-forwardsolve} we show the evolution of the deformation. That is, we solve the forward problem from $t=0$ to $t=1$ based on the found velocity field $v(x)$. We limit this visualization to the sagittal view. We display intermediate states for $t^j\in\{0,0.25,0.5,0.75,1\}$. In the top row, we show the mismatch between the state variable $m$ and the reference image $m_R$. In the two bottom rows, we show the state variable itself and the state variable with a grid in overlay to illustrate the deformation. The configuration at $t=1$ is equivalent to the results visualized in \figref{f:results-sphere-to-bowl}.

\ipoint{Observations}: The results reported in~\figref{f:results-sphere-to-bowl} and~\figref{f:results-sphere-to-bowl-forwardsolve} illustrate that we can recover large deformations. We reduce the mismatch between the transported template image and the reference image significantly. The mismatch is almost zero everywhere after registration. Only at the sharp corners of the bowl and along its boundary we can observe residual differences between the deformed template image and the reference image (dark areas). In~\cite{Mang:2017a} (in a 2D setting) we experimentally observed that using time-dependent velocities may help to further reduce this mismatch.

Careful visual inspection of the evolution of the deformation (see \figref{f:results-sphere-to-bowl-forwardsolve}) suggests that we obtain a well behaved map $y$. This is confirmed by the values for the determinant of the deformation gradient; the determinant of the deformation gradient does not change sign and is non-zero for all $x$. The extremal values are ($\max,\min)= (\num{5.302593401e-02},\num{2.732119508e+01})$. We conclude that the deformation map is locally diffeomorphic (up to numerical accuracy) as judged by these values.

\subsubsection{Smooth Synthetic Data}

We report results for a smooth synthetic test problem to demonstrate the performance of our solver under idealized conditions (smooth data and smooth velocity field without the presence of noise).

\ipoint{Setup}: The template image is given by $m_T(x) = \frac{1}{3}(\sin x_1 \sin x_1 + \sin x_2 \sin x_2 + \sin x_3 \sin x_3)$. The velocity field is given by $v(x) = 0.5(v_1(x),v_2(x),v_3(x))^\T$, $v_1(x) = \sin x_3 \cos x_2 \sin x_2$, $v_2(x)= \sin x_1 \cos x_3 \sin x_3$, $v_3(x) = \sin x_2 \cos x_1 \sin x_1$. The reference image $m_R(x)$ is computed by transporting the template image with $v$. We solve the inverse problem on a grid size of $256^3$ (\num{50331648} unknowns). We set the tolerance for the relative reduction of the gradient to $\num{1E-5}$ and the absolute tolerance for the gradient to $\num{1E-8}$ (not reached). We use an $H^2$ regularization model with a regularization parameter of $\num{1E-4}$. We set the maximal number of Newton iterations to 50 and the maximal number of Krylov iterations to 100.

\ipoint{Results}: We report the convergence behaviour of our solver in \figref{f:convergence-inverse-solve} (top row). We show (from left to right) the trend of the relative values (normalized with respect to the values obtained at the first iteration) for the mismatch, the $\ell^2$-norm of the gradient, and objective values with respect to the Gauss--Newton iteration index. The last graph shows the residual of the PCG per Gauss--Newton iteration (from Gauss--Newton iteration 1 (light blue) to Gauss--Newton iteration 5 (dark blue)), respectively.

\ipoint{Observations}: Our solver converges after six Gauss--Newton iterations. Overall, we require 73 Hessian matrix vector products and a total of 168 PDE solves. We require 1, 1, 4, 9, 23, and 35 PCG iterations per Gauss--Newton iteration, respectively (see also \figref{f:convergence-inverse-solve}; right column at the top). We reduce the gradient by $\|\di{g}_6\|_2/\|\di{g}_0\|_2 =  \num{6.548846578000E-06}$ to an absolute value of $\|\di{g}_6\|_2 =   \num{9.384663828245E-09}$ (the initial norm of the gradient is $\|\di{g}_0\|_2 =  \num{1.433025452111E-03}$). The runtime of our solver is \SI{4.171917e+02}{\second} ($\sim$\SI{6}{\minute} \SI{57}{\second}) on one node with 20 cores. We reduce the mismatch from \num{3.166924e-01} to \num{2.198840e-04}. Considering the plots in \figref{f:convergence-inverse-solve}, we can observe that we  mismatch, the gradient, and the objective drop rapidly. We converged to a stationary value of the objective functional and the mismatch after only 4 iteration. The relative reduction of the gradient at iteration 4 is $\|\di{g}_6\|_2/\|\di{g}_0\|_4= \num{5.032750773957E-03}$ ($\|\di{g}_4\|_2 = \num{7.212059953214E-06}$). During the last three iterations the gradient is still significantly reduced, but the mismatch and the objective value only change marginally. Considering the residual for the PCG iterations, we can see that the residual not only drops fast, but also that at each new Gauss--Newton iteration, we start almost at the same residual we terminated the former iteration with. That is, by using a superlinear forcing sequence we do not seem to oversolve for the Gauss--Newton step for this particular test problem. This changes when we move to real data.

\subsubsection{Real Brain Imaging Data}

We report results for real brain imaging data in a challenging multi-subject registration problem (i.e., the registration of brain imaging data of different individuals).

\begin{figure}
\includegraphics[width=\columnwidth]{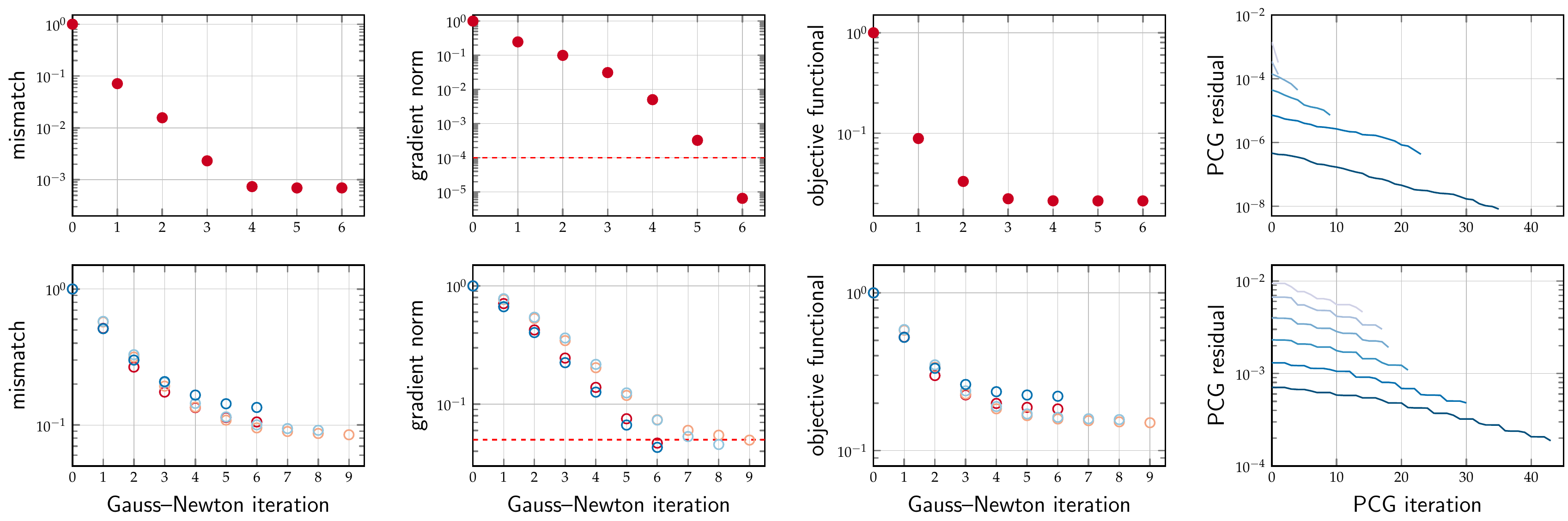}
\caption{Convergence for a synthetic problem (top row) and MRI brain data (bottom row; the different colors for the first three plots from the left correspond to different datasets). From left to right: We report the trend of ($i$) the relative mismatch, ($ii$) the relative gradient norm (the dashed red line indicates the tolerance for terminating the optimization), and ($iii$) the objective functional per Gauss--Newton iteration, respectively. We also show (right column) the trend of the residual for the PCG iterations (the different shades of blue represent the trend of the PCG residual for each Gauss--Newton iteration; the iteration count increases from bright to dark).}
\label{f:convergence-inverse-solve}
\end{figure}

\ipoint{Setup}: The images are taken from the \emph{non-rigid image registration evaluation project}~\cite{Christensen:2006a}. This repository contains 16 rigidly aligned T1-weighted MRI brain data sets of different individuals with a grid size of $256\times300\times256$ (\num{58982400} unknowns).\footnote{Additional information on the data sets, the imaging protocol, and the pre-processing can be found in~\cite{Christensen:2006a} and at \url{http://www.nirep.org/}.} We show one dataset of this repository in \figref{f:convergence-inverse-solve} (dataset na01 (reference image) and dataset na03 (template image)). We consider the datasets na01 through na05 for these experiments; na01 serves as the reference image for all experiments.

We consider an $H^1$-regularization model in combination with a penalty on the divergence of the velocity field ($H^1$-div regularization model). The regularization parameter for the Sobolev norm for the velocity is set to $\num{1E-2}$.  The regularization parameter for the penalty on the divergence of the velocity is set to $\num{1E-4}$. These values were chosen by experimentation. We set the maximal number of Newton iterations to 50 and the maximal number of Krylov iterations to 100. We use a tolerance of $\num{5E-2}$ for the relative change of the gradient and set the absolute tolerance for the gradient to $\num{1E-6}$ (not reached).

\begin{remark}
We use tighter tolerances since we found by experimentation that decreasing the gradient further, significantly increases the runtime. The reduction in mismatch associated with this increase in runtime does, from our perspective, not justify the use of smaller tolerances for the reduction of the gradient.
\end{remark}

\begin{figure}
\includegraphics[width=\columnwidth]{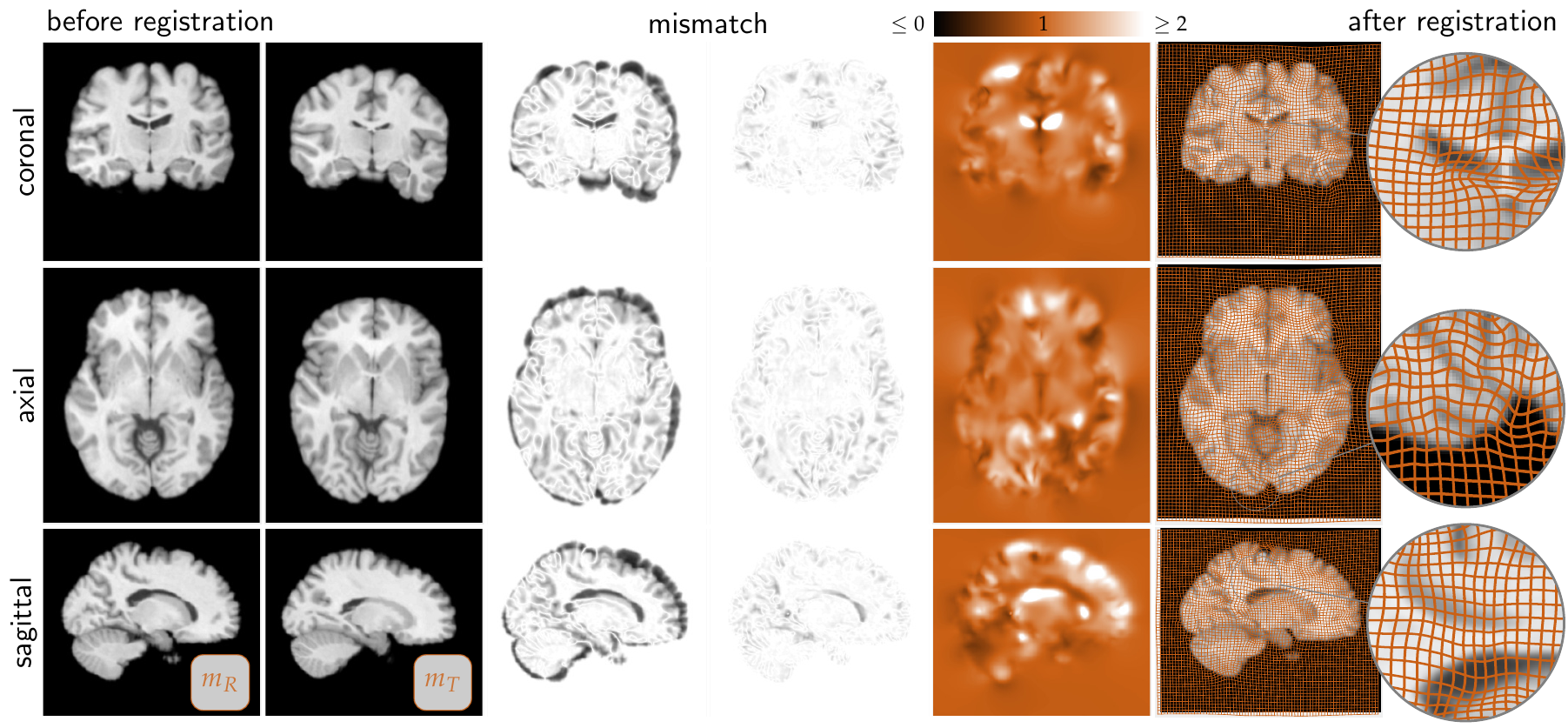}
\caption{Results for a real-world registration problem. We show (from left to right) the reference image $m_R$, the template image $m_T$, the residual differences before registration, the residual differences after registration, a pointwise map of the determinant of the deformation gradient, and the deformed template image with an illustration of the deformed grid in overlay. Each row corresponds to a different view of the 3D volume (from top to bottom: coronal, axial, and sagittal). Notice that the deformation map illustrated on the left is the pullback map (depicts, where points originate from), whereas the map for the determinant of the deformation gradient corresponds to the inverse of this deformation map.}
\label{f:nirep-3D-results-na03-to-na01-h1-ric}
\end{figure}

\ipoint{Results}: We report the convergence behaviour of our solver in \figref{f:convergence-inverse-solve} (top row). We show results for four registrations (na0$i$, $i=2,3,4,5$, to na01). We display (from left to right) the trend of the relative values (normalized with respect to the values obtained at the first iteration) for the mismatch, the $\ell^2$-norm of the gradient, and objective values, with respect to the Gauss--Newton iteration index. The last graph shows the residual of the PCG per Gauss--Newton iteration (from Gauss--Newton iteration 1 (light blue) to Gauss--Newton iteration 5 (dark blue)), respectively. The last plot is only for the registration of the dataset na02 to the dataset na01. We provide additional measures of performance in \tabref{t:nirep-results}. We illustrate results for an exemplary dataset (na03 to na01) in~\figref{f:nirep-3D-results-na03-to-na01-h1-ric}.

\begin{table}
\caption{Registration performance for MRI brain data of size $256\times300\times256$ (\num{58982400} unknowns). The results correspond to those reported in~\figref{f:convergence-inverse-solve}. We report the number of Gauss--Newton iterations $k$, the number of Hessian matvecs, the number of PDE solves, the relative mismatch after registration, the $\ell^2$-norm of the gradient, the relative reduction of the $\ell^2$-norm of the gradient, and the runtime. The results are obtained on a single node of CACDS's Opuntia cluster (2-socket Intel Xeon E5-2680v2 at \SI{2.8}{\giga\hertz} with 10 cores/socket and \SI{64}{\giga\byte} memory). We do not perform any grid, scale, or parameter continuation.\label{t:nirep-results}}
\centering
\begin{tabular}{llllllll}\toprule
dataset & $k$ & matvecs & PDE solves &                 mismatch &         $\|\di{g}_k\|_2$ & $\|\di{g}_k\|_2/\|\di{g}_0\|_2$ & runtime \\\midrule
na02    & 6   &     143 &        308 & \num{1.053063174936E-01} & \num{9.407527189929E-03} &        \num{4.683876639045E-02} & \SI{9.696632e+02}{\second} (\SI{16}{\minute}\phantom{00}\SI{9}{\second}) \\
na03    & 9   &     255 &        563 & \num{8.480884587284E-02} & \num{4.877870063395E-04} &        \num{4.962401062151E-02} & \SI{1.838740e+03}{\second} (\SI{30}{\minute}\phantom{0}\SI{38}{\second}) \\
na04    & 8   &     203 &        434 & \num{9.145664973708E-02} & \num{4.571084238661E-04} &        \num{4.556013285644E-02} & \SI{1.410249e+03}{\second} (\SI{23}{\minute}\phantom{0}\SI{30}{\second}) \\
na05    & 6   &     137 &        296 & \num{1.346620316391E-01} & \num{4.072076642777E-04} &        \num{4.308053573814E-02} & \SI{9.500624e+02}{\second} (\SI{15}{\minute}\phantom{0}\SI{50}{\second}) \\\bottomrule
\end{tabular}
\end{table}

\ipoint{Observations}: Our solver converges after 6 to 9 Gauss--Newton iterations for these datasets for a tolerance of $\num{5e-2}$ for the relative reduction of the gradient. We require between 296 and 563 PDE solves (this includes the solves for the evaluation of the objective functional, the gradient, and the Hessian matvecs) and 137 up to 255 Hessian matvecs (see \tabref{t:nirep-results}). The mismatch is reduced by a factor of \num{1.346620316391E-01} up to \num{8.480884587284E-02}. We can see that the mismatch drops quickly during the first few iterations. After about 6 iterations the reduction stagnates. The same is true for the objective functional. From these plots we can also see that further reducing the gradient will not lead to a significant reduction in the mismatch (which is what we care about from a practical perspective). For some of the runs we can also see that the reduction in the gradient starts to stagnate. If we use a smaller tolerance for the gradient, the runtime will increase significantly. Overall, we obtain a runtime of 15 minutes up to 30 minutes. This is not competitive with certain existing methods for diffeomorphic registration (\cite{Vercauteren:2009a} is a prominent example). We note, that these methods do, in general, not monitor the reduction of the gradient. They just perform a fixed number of iterations. Another key difference is, that these approaches perform a grid continuation, something we have not considered here. We expect that this will significantly reduce the runtime, given our solver has a complexity of $\mathcal{O}(n)$, where $n$ is the number of unknowns.

In \figref{f:convergence-inverse-solve} we can also see that for the considered dataset the residual for the PCG iterations drops steadily, but slowly (these results are for the registration of the dataset na02 to dataset na01). We require 14, 17, 18, 21, 30, and 43 PCG iterations per Gauss--Newton iteration, respectively. Qualitatively, we can observe that the template and reference image are in good agreement (see \figref{f:nirep-3D-results-na03-to-na01-h1-ric}). The deformation map is locally diffeomorphic as judged by the value for the determinant of the deformation gradient. The visualization of the deformed grid also illustrates a well behaved deformation map.

\subsection{Data Assimilation}
\label{s:experiments:tumor}

We study the performance of our formulation for synthetic test problems, for varying noise perturbations and detection thresholds.

\ipoint{Setup}: We apply different levels of noise and consider partial observations of the computed states. We assume observations at two time points ($t=0$ and $t=1$). We consider different thresholds for the partial observations since there is no consensus in the literature on how to model the imaging operator (i.e., how to relate the observed imaging abnormalities to the computed cell density)~\cite{Swanson:2008a,Harpold:2007a,Mang:2012a}. Using this synthetic test case allows us to study the inversion quality as a function of noise in the observation and as a function of the detection threshold. We use the atlas in~\cite{Cocosco:1997a} to match our implementation to virtual brain anatomy. The diffusion tensor imaging dataset is taken from~\cite{Mori:2008a}. We consider a multifocal tumor case. We invert for the initial condition and the diffusion coefficient.

\begin{figure}
\centering
\includegraphics[width=\columnwidth]{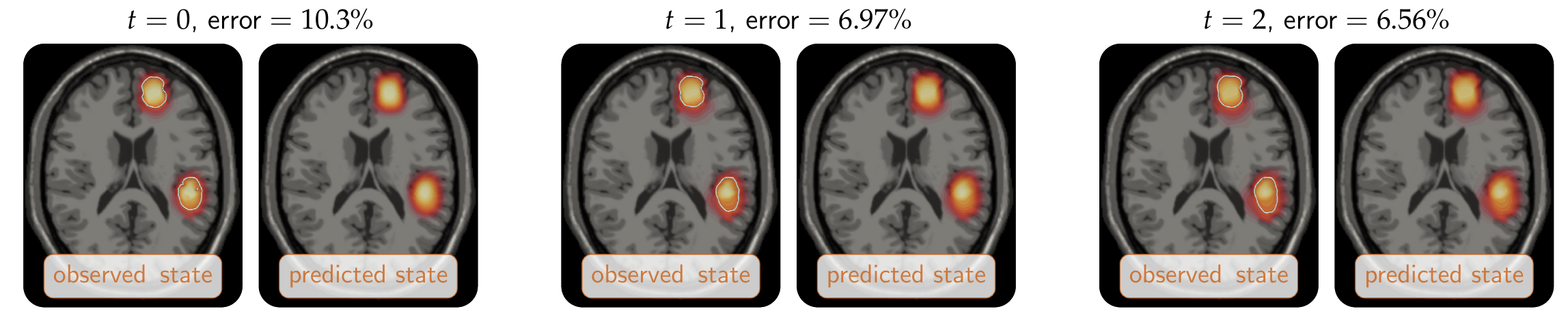}
\caption{Inversion results for a multifocal brain tumor (two dimensional problem). We show three observations for the time points $t=0,1,2$ (in dimensional form: $0$, $5$ and $10$ months; left image: `'observed state`') and the corresponding reconstructions (right image: `'predicted state`'). The noise level for the observations is set to $5\%$. We show isolines that correspond to the threshold we use to model the partial observation (white lines in the `'observed state`'; the threshold is $0.2$). The relative error between predicted and observed state is $10.3\%$, $6.97\%$, and $6.56\%$, respectively. Images modified from~\cite{Gholami:2016a}.}
\label{f:results-multifocal-tumor-inversion}
\end{figure}

\ipoint{Results}: We showcase representative results in~\figref{f:results-multifocal-tumor-inversion} (2D setting) and report relative reconstruction errors in~\tabref{t:tumor-inversion}. We report reconstruction errors for three timepoints, $t=0$, $t=1$, and $t=2$. The first two timepoints are used as observations in the inversion. The results presented in this study were originally reported in~\cite{Gholami:2016a}.

\begin{table}
\centering
\caption{Detailed reconstruction results for a multifocal tumor case (see \figref{f:results-multifocal-tumor-inversion}). We report relative errors for the computed diffusion coefficient and the predicted and observed tumor distribution at different points in time as a function of the noise and the detection threshold for the observed data. Results originally reported in~\cite{Gholami:2016a}.\label{t:tumor-inversion}}
\scriptsize
\begin{tabular}{rrrrrr}\toprule
          &       & \multicolumn{4}{l}{reconstruction error}                      \\\midrule
threshold & noise & coefficient   & $t=0$         & $t=1$         & $t=2$         \\\midrule
0.1       &  1\%  & \num{1.0e-02} & \num{9.7e-02} & \num{6.1e-02} & \num{5.5e-02} \\
          &  5\%  & \num{2.3e-02} & \num{1.0e-01} & \num{6.6e-02} & \num{6.1e-02} \\
          & 10\%  & \num{3.8e-02} & \num{1.1e-01} & \num{7.9e-02} & \num{7.4e-02} \\
0.2       &  1\%  & \num{8.0e-02} & \num{1.0e-01} & \num{7.2e-02} & \num{6.6e-02} \\
          &  5\%  & \num{8.0e-02} & \num{1.1e-01} & \num{7.8e-02} & \num{7.3e-02} \\
          & 10\%  & \num{8.1e-02} & \num{1.2e-01} & \num{9.2e-02} & \num{8.7e-02} \\
0.3       &  1\%  & \num{7.4e-02} & \num{1.2e-01} & \num{9.0e-02} & \num{8.5e-02} \\
          &  5\%  & \num{6.5e-02} & \num{1.2e-01} & \num{9.6e-02} & \num{9.1e-02} \\
          & 10\%  & \num{5.5e-02} & \num{1.4e-01} & \num{1.1e-01} & \num{1.1e-01} \\
0.4       &  1\%  & \num{1.7e-01} & \num{1.3e-01} & \num{1.0e-01} & \num{9.4e-02} \\
          &  5\%  & \num{1.7e-01} & \num{1.3e-01} & \num{1.0e-01} & \num{9.7e-02} \\
          & 10\%  & \num{1.7e-01} & \num{1.4e-01} & \num{1.1e-01} & \num{1.1e-01} \\
\bottomrule
\end{tabular}
\end{table}

\ipoint{Observations}: Qualitatively, we can see that the reconstruction (predicted state) is in good agreement with the observed sate (see \figref{f:results-multifocal-tumor-inversion}). The relative error between predicted and observed state for the images depicted in \figref{f:results-multifocal-tumor-inversion} is $10.3\%$, $6.97\%$, and $6.56\%$, respectively. The images shown in this figure correspond to the results for 5\% noise reported in \tabref{t:tumor-inversion}. Here, we report results for different time points, different noise levels, and different detection thresholds. The error increases as we increase the detection (observation) threshold. This is expected, since a higher detection threshold means that we observe less data. We can also see that the error increases as we increase the noise level. This also follows our intuitive understanding of the problem. The relative error for the diffusion coefficient ranges from \num{1.0e-02} for a detection threshold of 0.1 and 1\% noise to \num{1.7e-01} for a detection threshold of 0.4 and 10\% noise.

\section{Conclusions}
\label{s:conclusions}

Our research on the design of effective algorithms for model-driven data analytics relates to numerous areas in computational sciences and engineering beyond brain-tumor imaging, including biophysics inversion in medical imaging in general, computational geosciences, numerical weather prediction and atmospheric sciences. PDE-constrained optimization problems pose significant mathematical challenges. We have to deal with space-time, nonlinear, strongly coupled forward operators. The inverse problem is high-dimensional, nonlinear, nonconvex, and inherently ill-posed. This makes parallel algorithms and high-performance computing platforms necessary for computational tractability.

We have discussed recent advances in PDE-constrained optimization problems applied to medical image analysis. We considered three problems in medical image analysis---data assimilation in brain tumor imaging, motion estimation in cardiac imaging, and diffeomorphic image registration. We gave a comprehensive review of existing approaches in the field. We have showcased numerical results for real and synthetic data to demonstrate the performance and capabilities of our methods. The design of a distributed-memory Newton--Krylov solver for these type of problems is exclusive to our group~\cite{Mang:2015a,Mang:2017c,Gholami:2016a,Mang:2016c,Gholami:2017a}. We implement parallelism based on the message passing interface. Our solver is designed to scale on high-performance computing platforms~\cite{Mang:2016c,Gholami:2017a}. In general, we use the same numerical approach to solve the optimization problem for all formulations considered in the present manuscript. However, we note that we have to implement tailored numerics for each individual problem. That is, the form of the mixed-type PDE-operators that appear in our optimality systems dictate the design of the solvers for the individual subblocks.

In our past work, we have identified several components of our numerical schemes that can be improved. For example, for the diffeomorphic image registration problem we found that our preconditioners are not $\beta$-independent~\cite{Mang:2017c}. Designing an effective preconditioner for vanishing regularization parameters requires more work. This is something we are actively working on at the moment. For $H^1$-regularity for the velocity field, we observe good performance. However, if we increase the order of the regularization operator, the performance of our solver deteriorates.

We have deployed a solver that cannot only handle problems of unprecedented scale (for example, $\sim$200 billion unknowns for the diffeomorphic registration problem), but also solve clinically relevant problems in real-time (for example, under 2 seconds using 512 cores for the diffeomorphic registration problem)~\cite{Gholami:2017a}. As of today, this requires access to supercomputing platforms. To make these methods more useful for practitioners that may not have access to these kind of resources, we will work on improving our single node performance as well as port our solvers to GPU architectures.

Another key challenge for our future work are uncertainties in the data, the parameters, the algorithms, and the forward model. While predictive mathematical models of complex biological phenomena hold tremendous promise for revolutionizing computational medicine, their practical use is limited unless we rigorously quantify these uncertainties (i.e., provide measures of the degree of confidence with which our methods predict particular quantities of interest). The methods for data assimilation in brain tumor imaging discussed here have to be considered as preliminary. A significant amount of work remains ahead of us before this type of technology can be used in clinical decision making. The tumor model is in its present form not sophisticated enough to allow for tumor growth prediction. It neglects key elements of tumor progression. In future work, we need to investigate more complex multi-species models that, e.g., take into account mass effect, edema, necrosis, angiogenesis, chemotaxis, or haptotaxis~\cite{Garcke:2016a,Toma:2012e,Schuetz:2013a,HawkinsDaarud:2013a,Lima:2014a,Hogea:2007a,Hogea:2008b}. The integration of more complex models into our framework poses significant computational and theoretical challenges. Accurate and scalable schemes to evaluate and approximate adjoint, gradient, and Hessian operators are essential, not only in a deterministic setting, but also in the context of Bayesian inversion~\cite{Geweke:1999a,Geweke:2003a,Martin:2012a,Petra:2014a}.

In our most recent work, we have augmented our formulation for diffeomorphic registration with our formulation for biophysics inversion~\cite{Scheufele:2017a}; related approaches have been described in~\cite{Hogea:2008a,Gooya:2012a,Zacharaki:2009a}. Integrating these methods into a common mathematical framework for biophysics-based image analysis that combines variational methods with supervised machine learning for glioma segmentation~\cite{Mang:2017d} has shown to be a promising approach in the context of glioma brain tumor segmentation~\cite{Menze:2015a}.

\begin{appendix}

\section{Computing $\det \igrad y$}
\label{s:deformation-gradient}

We do not differentiate the deformation map $y$, but transport $\psi \defeq \det\igrad \phi$ where $y^{-1}(x) = \phi(x,t=1)$. That is, we solve
\begin{subequations}
\label{e:transport-jac}
\begin{align}
\p_t \psi + v\cdot\igrad\psi
& = \psi \idiv v
&& \text{in}\;\;\Omega\times(0,1]
\\
\psi &= 1
&& \text{in}\;\;\Omega\times\{0\}
\end{align}
\end{subequations}

\noindent with periodic boundary conditions on $\partial\Omega$.

\begin{remark}
If we start from an identity map $\operatorname{id}$ (zero velocity field) the determinant of the deformation gradient will have a value of one. Hence, we require that the determinant of the deformation gradient remains strictly positive for all $x \in \Omega$ for the computed deformation map to be locally diffeomorphic. Strictly speaking, we also require that the map $y$ is smooth, and has a smooth inverse for $y$ to be a diffeomorphism. If we use (slightly more than) $H^2$-regularity for $v$ we meet these requirements from a theoretical perspective (assuming that the images are smooth)~\cite{Beg:2005a,Vialard:2012a}. However, we note that ensuring that the discrete deformation map $y$ (which does not appear in our formulation) is a diffeomorphism, requires sophisticated discretization schemes (see, e.g.,~\cite{Burger:2013a}). In our formulation, we transport the data. Since the deformation map $y$ does not appear explicitly, we decided to compute the determinant of the deformation gradient by computing the solution of the transport problem in \secref{s:deformation-gradient}. We found that this is more stable than computing $\det \igrad y$ from $y$ using our numerical scheme.
\end{remark}

\section{Gauss--Newton Krylov Algorithm}
\label{s:gn-krylov-method}

Here, we showcase the individual steps necessary for solving~\eqref{e:pde-opt} in~\algref{a:outer-iteration} (Newton iterations; \emph{outer iterations}) and~\algref{a:inner-iteration} (solution of the system in~\eqref{e:reduced-space-kkt-system}, \emph{inner iterations}). We refer to~\cite{Mang:2015a,Mang:2017c,Gholami:2016a} for details on the particular implementation of our nonlinear solver. We rely on PETSc's TAO package~\cite{Balay:2016a,Munson:2017a} for our distributed-memory solver~\cite{Mang:2016c,Gholami:2017a} (see \secref{s:implementation-details} for details). The steps are the same as outlined in~\algref{a:outer-iteration} and~\algref{a:inner-iteration}, respectively.

\begin{algorithm}
\caption{Inexact Newton--Krylov method.}
\label{a:outer-iteration}
\algadjust
\begin{algorithmic}[1]
\STATE{$\di{w}_0 \leftarrow \digr{0}$, $k\leftarrow0$}
\STATE{$\di{m}_0 \leftarrow$ solve state equation forward in time given $\di{w}_0$}
\STATE{$\digr{\lambda}_0 \leftarrow$ solve adjoint equation backward in time given $\di{w}_0$ and $\di{m}_0$}
\STATE{evaluate $\F{J}^h$ and $\di{g}_0$ given $\di{m}_0$, $\digr{\lambda}_0$ and $\di{w}_0$}
\WHILE{$\|\di{g}_k\|^2_2 > \|\di{g}_0\|_2^2\epsilon_{\text{opt}}$}
    \STATE{$\di{\tilde{w}}_k\leftarrow$ solve $\di{H}_k \di{\tilde{w}}_k = -\di{g}_k$ given $\di{m}_k$, $\digr{\lambda}_k$, $\di{w}_k$, and $\di{g}_k$}
    \STATE{$\alpha_k \leftarrow$ perform line search on $\di{\tilde{w}}_k$ subject to Armijo condition}
    \STATE{$\di{w}_{k+1} \leftarrow \di{w}_k + \alpha_k\di{\tilde{w}}_k$}
    \STATE{$\di{m}_{k+1} \leftarrow$ solve state equation forward in time given $\di{w}_{k+1}$}
    \STATE{$\digr{\lambda}_{k+1} \leftarrow$ solve adjoint equation backward in time given $\di{w}_{k+1}$ and $\di{m}_{k+1}$}
    \STATE{evaluate $\F{J}^h$ and $\di{g}_{k+1}$ given $\di{m}_{k+1}$, $\digr{\lambda}_{k+1}$ and $\di{w}_{k+1}$}
    \STATE{$k \leftarrow k + 1$}
\ENDWHILE
\end{algorithmic}
\end{algorithm}

\begin{algorithm}
\caption{Newton step. We illustrate the solution of the reduced KKT system~\eqref{e:reduced-space-kkt-system} using a PCG method at a given outer iteration $k\in\ns{N}$.}
\label{a:inner-iteration}
\algadjust
\begin{algorithmic}[1]
\STATE{input: $\di{m}_k$, $\digr{\lambda}_k$, $\di{w}_k$, $\di{g}_k$, $\di{g}_0$}
\STATE{set $\eta_k \leftarrow\min(0.5,(\|\di{g}_k\|_2/\|\di{g}_0\|_2)^{1/2})$, $\di{\tilde{w}}_0 \leftarrow \digr{0}$, $\di{r}_0 \leftarrow - \di{g}_k$}
\STATE{$\di{z}_0 \leftarrow$ apply preconditioner to $\di{r}_0$}
\STATE{$\di{s}_0 \leftarrow \di{z}_0$, $\iota\leftarrow 0$}
\WHILE{$\iota < n$}
    \STATE{$\di{\tilde{m}}_\iota \leftarrow$ solve inc. state eq. forward in time given $\di{m}_k$, $\di{w}_k$ and $\di{\tilde{w}}_\iota$}
    \STATE{$\digr{\tilde{\lambda}}_\iota \leftarrow$ solve inc. adjoint eq. backward in time given $\digr{\lambda}_k$, $\di{w}_k$, $\di{\tilde{m}}_\iota$ and $\di{\tilde{w}}_\iota$}
    \STATE{$\di{\tilde{s}}_\iota \leftarrow$ apply $\di{H}_\iota$ to $\vect{s}_\iota$ given $\digr{\lambda}_k$, $\di{m}_k$, $\di{\tilde{m}}_\iota$ and $\digr{\tilde{\lambda}}_\iota$ (Hessian matvec)}
    \STATE{$\kappa_\iota \leftarrow\langle\di{r}_\iota,\di{z}_\iota\rangle/\langle\di{s}_\iota,\di{\tilde{s}}_\iota\rangle$,
            \quad $\di{\tilde{w}}_{\iota+1} \leftarrow \di{\tilde{w}}_\iota + \kappa_\iota\di{s}_\iota$,
            \quad $\di{r}_{\iota+1}\leftarrow\di{r}_\iota-\kappa_\iota\di{\tilde{s}}_\iota$}
    \STATE{\textbf{if} $\|\di{r}_{\iota+1}\|_2 < \eta_k$ \textbf{break}}
    \STATE{$\di{z}_{\iota+1} \leftarrow$ apply preconditioner to $\di{r}_{\iota+1}$}
    \STATE{$\mu_\iota \leftarrow \langle\di{z}_{\iota+1},\di{r}_{\iota+1}\rangle/\langle\di{z}_\iota,\di{r}_\iota\rangle$, \quad $\di{s}_{\iota+1} \leftarrow \di{z}_{\iota+1} + \mu_\iota\di{s}_\iota$, \quad $\iota \leftarrow \iota + 1$}
\ENDWHILE
\STATE{output: $\di{\tilde{w}}_k \leftarrow\di{\tilde{w}}_{\iota+1}$}
\end{algorithmic}
\end{algorithm}

\section{Newton Step}
\label{s:newton-step}

Here, we describe the Newton step for the two problems described in~\secref{s:formulation-reg} and~\secref{s:formulation-tumor}. The Newton step is obtained by computing the second variation of~\eqref{e:lagrangian-reg}. We will see that the structure of the PDE operators that appear in the \emph{reduced space} Hessian is very similar to the first-order optimality systems in~\secref{s:optimality-systems:reg}.

\subsection{Diffeomorphic Registration}
\label{s:newton-step-reg}

The nonlinearity, ill-posedness, and the dimensionality of the search space makes the solution of this variational problem computational challenging. Most available implementations use first-order iterative schemes to compute a minimizer to the variational optimization problem (see, e.g.,~\cite{Avants:2011a,Beg:2005a,Chen:2011a,Ha:2010a,Hart:2009a,Vialard:2012a}). To increase the rate of convergence, they typically do not use $g_v$ in their iterative scheme, but the gradient in the Sobolev space induced by the regularization operator $\D{A}$, i.e.,
\begin{align*}
\tilde{g}_v = v + (\beta\D{A})^{-1}\int_0^1\lambda\igrad m\d{t}.
\end{align*}

First-order methods are usually inefficient; they require a lot of iterations to converge to a specific tolerance; we have analyzed this in~\cite{Mang:2015a}. An alternative is to apply variants of Newton's method to the KKT system. To derive the operators for applying Newton's method we have to compute the second variations of~\eqref{e:lagrangian-reg}. We arrive at the incremental state and adjoint equation
\begin{subequations}
\label{e:newton-step}
\begin{align}
\p_t \tilde{m} + \igrad \tilde{m} \cdot v
 + \igrad m \cdot \tilde{v} & = 0
&&{\rm in}\;\; \Omega \times (0,1],
\label{e:inc-state}
\\
\tilde{m} &= 0
&&{\rm in}\;\; \Omega \times\{0\},
\label{e:inc-init-state}
\\
-\p_t \tilde{\lambda} - \idiv (v\tilde{\lambda}
+ \tilde{v}\lambda) & = 0
&&\text{in}\;\; \Omega \times [0,1),
\label{e:inc-adjoint}
\\
\tilde{\lambda} &= - \tilde{m}
&&{\rm in}\;\; \Omega \times\{1\},
\label{e:inc-final-adj}
\\
\intertext{and the expression for the reduced space Hessian matvec}
\D{H}(v)\tilde{v} \defeq \beta \D{A}[\tilde{v}]
+ \D{K}[\int_0^1\tilde{\lambda}\igrad m + \lambda\igrad \tilde{m}\d{t}\;]
&&&\text{in}\;\; \Omega,
\label{e:inc-control}
\end{align}
\end{subequations}

\noindent with periodic boundary conditions on $\p\Omega$. The incremental control variable $\tilde{v}$ corresponds to the search direction of our scheme. The transport equations for the incremental state and adjoint fields $\tilde{m}$ and $\tilde{\lambda}$ are hidden in the integro-differential operator in~\eqref{e:inc-control}; we have to solve~\eqref{e:inc-state} and~\eqref{e:inc-adjoint} every time we apply $\D{H}$ to a new vector $\tilde{v}$. We use a Gauss--Newton approximation to the true Hessian to ensure that the operator is positive definite far away from the optimum. This corresponds to dropping $\lambda$ in~\eqref{e:inc-adjoint} and~\eqref{e:inc-control}.

\subsection{Data Assimilation}
\label{s:newton-step-tumor}

The expressions for evaluating the Hessian operator are derived by computing the second variations for~\eqref{e:lagrangian-tumor}. The Hessian matvec is given by $\D{H}\tilde{p} \defeq \beta \tilde{p} - \Phi^\T \tilde{\lambda}_0$, $\D{H} : \ns{R}^{n_p} \rightarrow \ns{R}^{n_p}$. To be able to evaluate this expression we have to solve the incremental state and adjoint equations given by
\begin{subequations}
\begin{align}
\label{e:tumor:inc-state-pde}
\p_t \tilde{m} - \idiv k \igrad \tilde{m} - \rho(1-2m)\tilde{m}
& = 0
&& {\rm in}\;\Omega_B\times(0,1],
\\
\label{e:tumor:inc-state-initial}
\tilde{m}&=\Phi\tilde{p}
&& {\rm in}\;\Omega_B\times\{0\},
\\
\label{e:tumor:inc-adj-pde}
-\p_t \tilde{\lambda} - \idiv k \igrad \tilde{\lambda} - \rho(1-2m - 2\lambda)\tilde{\lambda}
& = 0
&& {\rm in}\;\Omega_B\times[0,1),
\\
\label{e:tumor:inc-adj-initial}
\tilde{\lambda}&= -\D{Q}^\T\D{Q} \tilde{m}
&& {\rm in}\;\Omega_B\times\{1\},
\end{align}
\end{subequations}

\noindent with Neumann boundary conditions on $\p\Omega_B$. For the Gauss--Newton approximation to the true Hessian $\lambda$ in \eqref{e:tumor:inc-adj-pde} needs to be dropped.

\end{appendix}

\end{document}